\documentclass[12pt]{amsart} 
\usepackage{amssymb}
\usepackage[french]{babel}

\newcommand{\got}{\mathfrak}

\newcommand{\lra}{\longrightarrow}
\newcommand{\dis}{\displaystyle}

\newcommand{\var}{\varepsilon}

\def\ind{\mathrm Ind}
\def\carre{{\vrule height5pt width5pt depth0pt}}
\def\qed{{\hfill\carre\vskip1em}}

\newcommand{\affiliationone}[1]{\noindent#1}
\newcommand{\classno}[1]{\medskip 2000 Mathematics Subject
  Classification~~ #1}

\title[Orbites Nilpotentes Sph\'eriques]
{Orbites Nilpotentes
  Sph\'eriques et Repr\'esentations unipotentes associ\'ees : Le cas 
$\bf SL_n$.}

\author{Herv\'e Sabourin}

\begin{document}

\maketitle

\begin{abstract} 
Let $G$ be a real simple Lie group, $\got g$
its Lie algebra. Given a nilpotent adjoint $G$-orbit $O$, the question
is to determine the irreducible unitary representations of $G$ that we can
associate to $O$, according to the orbit method. P.Torasso, in \cite{TO2}, gave a method to solve this problem
if $O$ is minimal. In
this paper, we study the case where $O$ is any spherical nilpotent
orbit of $sl_n({\mathbb  R})$, we construct, from $O$, a family of
representations of the two-sheeted covering of $SL_n({\mathbb R})$
with Torasso's method and, finally, we show that all these
representations are associated to the corresponding orbit.
\end{abstract}

\classno{20G05, 22E46, 22E47}

\section*{\bf 0. Introduction.}
 
Soit $G$ un groupe r\'eel simple, connexe et
simplement connexe,agissant sur son alg\`ebre de Lie $\got g$,via
l'action co-adjointe. ''La m\'ethode des orbites'', initi\'ee par
A.Kirillov dans le cadre des groupes nilpotents, a pour but d'essayer
de d\'ecrire le dual unitaire de $G$ \`a l'aide des orbites
coadjointes. Le probl\`eme est, donc, de d\'eterminer quelles
repr\'esentations unitaires irr\'eductibles  de $G$ on peut associer,
selon un sens \`a d\'efinir, \`a une $G$-orbite donn\'ee et on peut
s'int\'eresser, en particulier,  au cas o\`u cette orbite est nilpotente. 

 Il existe une mani\`ere naturelle d'associer une
repr\'esentation unitaire irr\'eductible d'un groupe simple $G$ \`a une
$G$-orbite nilpotente coadjointe $O$.
Soit, en effet, ${\got g}_{\mathbb
C}$ la complexifi\'ee de $\got g$ et $O_{\mathbb C}$ une $G_{\mathbb
C}$-orbite nilpotente. On suppose que 
l'intersection $O_{\mathbb C} \cap {\got g}$ est non vide. Dans ce cas,
cette intersection est une r\'eunion finie de $G$-orbites et on
suppose que $O$ est l'une de ces $G$-orbites. On sait, par ailleurs,
selon un r\'esultat obtenu ind\'ependamment par A.Joseph, W.Borho et
J.L.Brylinski (\cite{JO}, \cite{BB}), que tout id\'eal primitif de l'alg\`ebre enveloppante
$U({\got g})$ de ${\got g}_{\mathbb C}$ a pour vari\'et\'e des z\'eros 
l'adh\'erence  de
Zariski d'une et une seule orbite nilpotente complexe.
On dit alors qu'un \'el\'ement $\pi$ de $\widehat G$
est ``associ\'e'' \`a $O$ si la vari\'et\'e des z\'eros de  l'annulateur infinit\'esimal $Ann \pi$ de
$\pi$ dans $U({\got g})$ est l'adh\'erence de Zariski de
$O_{\mathbb C}$.
 
On peut remarquer imm\'ediatement que, si $\pi$ est
associ\'e \`a $O$ et si $GKdim (U({\got g}) / Ann\pi)$ d\'esigne la
dimension de Gelfand-Kirillov de l'alg\`ebre $(U({\got g}) / Ann\pi)$ , alors $\pi$ v\'erifie la condition $(GK)$
suivante~:
$$
GKdim (U({\got g}) / Ann\pi) = \dim O \leqno{({\bf GK})}
$$

Une repr\'esentation $\pi$ satisfaisant \`a $(GK)$ sera dite
``$GK$-associ\'ee'' \`a $O$. Par contre, il n'y a pas de raison a priori  qu'une
repr\'esentation ``$GK$-associ\'ee'' \`a $O$ soit ``associ\'ee'' \`a
$O$.
  
L'existence d'une repr\'esentation associ\'ee a \'et\'e \'etablie, dans
le cas o\`u l'orbite co-adjointe est nilpotente minimale, par
P.Torasso~\cite{TO2}, lorsque $G$ est de rang r\'eel sup\'erieur ou \'egal \`a
$3$, et par R.Brylinski et B.Kostant~\cite{BK} dans le cas g\'en\'eral, selon
des m\'ethodes totalement diff\'erentes. Dans le premier cas,
P.Torasso a mis au point un proc\'ed\'e permettant de donner une
construction explicite de le repr\'esentation, s'appuyant fortement
sur la notion d'orbite ``admissible'' au sens de M.Duflo~\cite{DU} et sur une
param\'etrisation du dual unitaire d'un groupe presque alg\'ebrique du
m\^eme M.Duflo. Dans le deuxi\`eme cas, R.Brylinski et B.Kostant ont
donn\'e une description du module de Harish-Chandra correspondant \`a
l'aide d'un proc\'ed\'e de quantification d'orbite.

Le point de vue qui nous int\'eresse, dans ce travail, est celui
de P.Torasso et M.Duflo. La m\'ethode utilis\'ee, sur laquelle nous
reviendrons en d\'etails plus tard, consiste en fait \`a
s'int\'eresser aux projections de l'orbite minimale sur les
sous-alg\`ebres paraboliques maximales standards, \`a construire une
repr\'esentation de chaque sous-groupe parabolique maximal associ\'e
\`a cette projection et \`a chaque donn\'ee d'admissibilit\'e et,
enfin, \`a consid\'erer le produit amalgam\'e de ces repr\'esentations
selon un  r\'esultat de J.Tits~\cite{SE}- C'est \`a cette occasion que l'on
a besoin de l'hypoth\`ese sur le rang- L'une des propri\'et\'es
essentielles servant \`a cette construction est le fait que l'orbite
minimale poss\`ede une $B$-orbite ouverte, $B$ \'etant un sous-groupe de
Borel de $G$.  Il semble donc naturel de s'int\'eresser aux $G$-orbites
nilpotentes r\'eelles dont la complexifi\'ee poss\`ede une $B_{\mathbb
C}$-orbite ouverte,
c'est-\`a-dire aux 
$G$-orbites nilpotentes sph\'eriques.

Dans ce travail nous nous int\'eressons aux orbites nilpotentes
sph\'eriques non minimales de $sl_n({\mathbb R})$, pour $n \geq
4$. On peut citer \`a ce propos Y.Flicker qui, dans \cite{FL}, a  \'etudi\'e
cette situation pour certaines de ces orbites.

Apr\`es avoir rappel\'e, dans les paragraphes {\bf 1} et {\bf 2},
les termes essentiels de la m\'ethode utilis\'ee, nous donnons une
description pr\'ecise des orbites consid\'er\'ees dans le paragraphe {\bf 3}. En
consid\'erant ensuite les projections sur les  paraboliques maximaux
d'une telle orbite $O$ et \`a l'aide d'un raisonnement par
r\'ecurrence on construit, \`a partir de $O$,  une famille de
repr\'esentations unitaires irr\'eductibles du rev\^etement \`a deux
feuillets de $SL_n({\mathbb R})$. Ceci fait l'objet des paragraphes
{\bf 4} et  {\bf 5}. 

Dans le paragraphe {\bf 6}, nous montrons, en utilisant une
r\'ealisation explicite de ces repr\'esentations, que celles-ci sont
toutes ``$GK$-associ\'ees'' \`a l'orbite $O$.

Enfin, dans le paragraphe {\bf 7}, nous montrons que ces
repr\'esentations sont aussi ``associ\'ees'' \`a $O$.

Je tiens \`a remercier tout particuli\`erement J.Y Charbonnel et
D.Vogan dont les  suggestions d\'eterminantes ont permis de
d\'emontrer le th\'eor\`eme 7.1., r\'esultat principal du paragraphe
{\bf 7}. 

\section{\bf Les param\'etrisations de Duflo.}

Nous allons rappeler dans ce
paragraphe deux param\'etrisations de M.Duflo, essentielles 
\`a la m\'ethode que nous allons
utiliser pour construire les repr\'esentations souhait\'ees. 

\subsection{} La premi\`ere param\'etrisation est celle
du dual unitaire d'un groupe presque alg\'e\-bri\-que r\'eel $P$ , d'alg\`ebre de Lie
$\got p$. Soit $q$ un \'el\'ement de ${\got p}^*$, ${\got p}(q)$ et $P(q)$,
respectivement le stabilisateur de $q$ dans ${\got p}$ et $P$. Soit
$B_q$ la forme bilin\'eaire altern\'ee sur $\got p$, d\'efinie par:
$\forall X,Y \in {\got p}, B_q(X,Y) = q([X,Y])$.
Dans \cite{DU}, Duflo introduit les notions suivantes:

\smallskip 

{\bf Definition 1.1.} 
{\it Une sous-alg\`ebre $\got b$ de $\got p$ est
dite de type fortement unipotent relativement \`a $q$ si $\got b$ est
alg\'ebrique, coisotrope rela\-tivement \`a la forme $B_q$, et si l'on
a: ${\got b} = {\got p}(q) + \ ^u{\got b}$.}

\smallskip

{\bf Definition 1.2.} 
{\it La forme $q$ est dite de type unipotent si les deux conditions
suivantes sont r\'ealis\'ees:

-Il existe un facteur r\'eductif de ${\got p}(q)$ contenu dans
$\ker q$.

-Il existe une sous-alg\`ebre de type fortement unipotent
relativement \`a $q$.}

Soit $R(q)$
un facteur r\'eductif de $P(q), {\got r}(q)$ son alg\`ebre de Lie.
L'espace ${\got p} / 
{\got p}(q)$ est muni d'une structure symplectique $R(q)$-
invariante, permettant de d\'efinir l'extension m\'etaplectique 
$R(q)^{\got p}$. Le noyau de cette extension admet 
un et un seul \'el\'ement non trivial, not\'e $e$. Selon G.Lion~\cite{LI-VE}, 
l'extension $R(q)^{\got p}$ se d\'ecrit de la mani\`ere suivante~:  
On se donne un Lagrangien $L$ de l'espace symplectique ${\got p} / {\got p}
(q)$ et on choisit une orientation $\widetilde{L}$ de $L$. A tout $x$ de 
$R(q)$, on associe l'orientation relative des Lagrangiens orient\'es  
$\widetilde{L}$ et $x.\widetilde{L}$ que l'on note 
$e(\widetilde{L},x.\widetilde{L})$. On pose, ensuite :
$$t(x)^2 =  e(\widetilde{L},x.\widetilde{L})$$
On obtient, alors :
$$R(q)^{\got p} = \{(x,t(x)), x \in R(q)\}$$

Soit $Y(q) = \{ \tau \in
\widehat{R(q)^{\got p}} \mid \tau(e)
= - Id \}$ et ${\mathbb E} = \{ (q ,\tau ) \mid q$ de type unipotent, $\tau \in
Y(q) \}$.
Le groupe $P$ op\`ere dans $\mathbb E$ et Duflo \'etablit, dans \cite{DU}, une
bijection de ${\mathbb E} / P$ sur $\widehat P$. L'image par cette bijection
d'un couple $(q,\tau)$ sera not\'ee $\pi_{q,\tau}$ et sera appel\'ee
$P$-repr\'esentation {\it "de type Duflo"}.

Consid\'erons, d'autre part, l'ensemble suivant :
$$Adm_P(q) = \{ \tau \in Y(q) / d\tau \ \hbox{est un multiple de} \ iq_{\mid 
{\got r}(q)} \}$$

\smallskip

{\bf D\'efinition 1.3 :} {\it
Si $Adm_P(q) \not= \emptyset$, l'orbite $P.q$ est dite admissible et 
l'ensemble $Adm_P(q)$ est l'ensemble des param\`etres d'admissibilit\'e 
de l'orbite.}

Soit donc $(q,\tau)$ un {\'e}l{\'e}ment de ${\mathbb E}$. $q$ \'etant
de type unipotent,
il existe des sous-alg{\`e}bres
de type fortement unipotent relativement {\`a} $q$,
$P(q)$-invariantes; choisissons-en une, soit $\got b$.
Posons ${}^u{\got b} = {\got a}$, le radical unipotent de $\got b$. Soit
$A$ le
sous-groupe analytique de $P$ correspondant;
alors $B = P(q). A = R(q) \times A$ est un sous-groupe ferm{\'e} de
$P$, d'alg{\`e}bre de Lie $\got b$. Soit $\mu$ la restriction de $q$ {\`a}
$\got a$.
$R(q)$ op{\`e}re dans $\got a$ en laissant stable $\mu$, de sorte
que l'extension
$R(q)^{\got a}$ est bien d{\'e}finie et que l'on peut associer \`a $\tau$
 un \'el\'ement $\tilde
\tau$ de $\widehat{R(q)^{\got a}}$, d\'efini par la formule suivante:
\begin{equation}
\forall (x, t'(x)) \in
R(q)^{\got a}, \ \tilde{\tau}(x,t'(x)) = {t(x) \over t' 
(x)}
\tau (x,t(x)) 
\end{equation}
(Cette formule ne d\'epend pas du choix du repr\'esentant $(x,t(x))$ 
de
$x$ dans $R(q)^{\got p}$).

Soit $T_\mu$ la classe de repr{\'e}sentations de $A$ associ{\'e}e {\`a} $
\mu$
par la correspondance de Kirillov, d'espace ${\got L}_\mu$, et $S_\mu$ la
repr{\'e}sentation m{\'e}taplectique correspondante.
On d{\'e}finit une repr{\'e}sentation du groupe $B$, not{\'e}e $\tau \otimes
S_\mu T_\mu$, dans le produit tensoriel de l'espace de $\tau$ et de l'espace
de $T_\mu$, soit $V_\tau \otimes {\got L}_\mu$, en posant~:
\begin{equation}
\forall x \in R(q),\ \forall y \in A, (\tau \otimes S_\mu
T_\mu) (xy) = {\tilde \tau}(x,t(x)) \otimes S_\mu
(x,t(x)).T_\mu(y). 
\end{equation}

Posons: 
\begin{equation}
\pi_{q,\tau, {\got b}} = \ind^P_B (\tau \otimes S_\mu T_\mu)
.
\end{equation}

D'apr\`es \cite{DU}, 3.16, on sait que si $\got b$ et $\got b'$ sont deux
sous-alg\`ebres de type fortement unipotent relativement \`a $q$,
$P(q)$-invariantes, alors les repr\'esentations $\pi_{q,\tau,{\got b}}$ et
$\pi_{q,\tau,{\got b}'}$ sont irr\'eductibles et \'equivalentes. La classe
d'\'equivalence de ces repr\'esentations est la $P$-repr\'esentation {\it
de type Duflo} $\pi_{q,\tau}$.

\subsection{} La deuxi\`eme param\'etrisation est celle de la
repr\'esentation coadjointe due encore \`a Duflo (\cite{DU}, chapitre
1). En reprenant les notations pr\'ec\'edentes, consid\'erons une
forme de type unipotent $q$ sur $\got
p$  et soit ${\got p}(q) = {\got r}(q) \oplus
\ ^u{\got p}(q)$. On introduit les ensembles suivants:
$$
\begin{array}{rl}
{\mathcal L}(q) &= \{\lambda \in {\got p}(q)^* \mid \lambda_{\mid
^u{\got
p}(q)} = q_{\mid ^u{\got p}(q)}\}, \\
{\mathcal D} &= \{ (q,\lambda) \mid \ q \ {\rm de \ type \ unipotent,} \ \lambda
\in {\mathcal L}(q) \}
\end{array}
$$
Notons que l'op\'eration ``restriction'' induit un isomorphisme de
${\mathcal L}(q)$ sur ${\got r}(q)^*$ et que le groupe $P$ op\`ere
naturellement sur $\mathcal D$.

Soit maintenant $(q,\lambda) \in {\mathcal D}$ et $\got b$ une
sous-alg\`ebre de type fortement unipotent relativement \`a $q$. Soit
$f \in {\got p}^*$ telle que:
$$f_{\mid ^u{\got b}} = q_{\mid ^u{\got b}}, \ \ f_{\mid {\got
p}(q)} = \lambda_{\mid {\got p}(q)}.$$
Duflo \'etablit les r\'esultats suivants:

- La $P$-orbite $P.f$ ne d\'epend pas des choix de $\got b$ et
$f$. On notera dor\'enavant $O_{q,\lambda}$ une telle orbite.

- L'application $(q,\lambda) \lra O_{q,\lambda}$ induit une
bijection de ${\mathcal D}/P$ sur ${\got p}^* /P$.

On remarquera, en particulier, que si $f$ est de type unipotent, 
alors, $P.f = O_{f,0}$. 

Rappelons deux  r\'esultats de Duflo utiles pour la
suite. 

\subsection{} Soit $U$ un sous-groupe unipotent de $P$ d'alg\`ebre de Lie
$\got u$, soit $u \in {\got u}^*, H = P(u), {\got v}$ une sous-alg\`ebre de
${\got u}, H$-invariante et co-isotrope relativement \`a $u$, $V$ le
sous-groupe correspondant. Soit
$v$ la restriction de $u$ \`a $\got v$. Les extensions
$H^{\got u}$ et $H^{\got v}$ 
 sont bien d\'efinies.  

Soit $\chi_u$ le caract\`ere de $U$ d\'efini par la forme
$u$. Soit $\tau$ une repr\'esentation de $H^{\got u}$ telle que :
$$\tau_{\mid U(u)} = \chi_u, \ \ \tau(e) = -Id$$
On associe \`a $\tau$ la  repr\'esentation $\widetilde{\tau}$ de 
$H^{\got v}$ d\'efinie par la  formule (1).
On a, alors, le r\'esultat suivant \cite{DU}, Lemme 17 :

\smallskip

{\bf  Proposition 1.1 :} {\it On a l'\'equivalence de
repr\'esentations suivante :
$$\ind_{P(u)V}^{P(u)U} ({\widetilde{\tau}} \otimes S_vT_v) 
  \simeq \tau
\otimes S_uT_u$$}

\subsection{}Les notations sont celles de 1.3. Soit $Q$ un
sous-groupe ferm\'e de $P$, contenant $U$. Soit $P_1 = P(u)^{\got u},
\ Q_1 = Q(u)^{\got u}, R_1$ une repr\'esentation unitaire de $Q_1$. On
pose :
$$T_1 = \ind_{Q_1}^{P_1}R_1, T'_1 = \ind_{Q(u)U}^{P(u)U} (R_1 \otimes S_uT_u)$$
M.Duflo, dans le chapitre II de \cite{DU},
 d\'emonstration de la
proposition II.15, prouve le r\'esultat suivant :

\smallskip

{\bf Proposition 1.2 :} {\it Les repr\'esentations $T'_1$ et $T_1
\otimes S_uT_u$ sont deux repr\'esentations \'equivalentes du groupe $P(u)U$.}

\section{\bf Syst\`emes de Tits et amalgames.} 
Nous allons rappeler,
dans ce paragraphe, les propri\'et\'es d'amalgame relatives aux
syst\`emes de Tits.

\smallskip

{\bf D\'efinition 2.1 :} {\it Un syst\`eme de Tits est un
quadruplet $(G,B,M',S)$ o\`u $G$ est un groupe, $B,M'$ deux
sous-groupes de $G$, $S$ une partie de $W = M' / B \cap M'$,
satisfaisant aux axiomes suivants~:

$(A_1)$  L'ensemble $B \cup M'$ engendre $G$ et $B \cap M'$
est un sous-groupe distingu\'e de $M'$.

$(A_2)$  L'ensemble $S$ engendre $W$ et se compose d'\'el\'ements
d'ordre $2$.

Pour tout $w \in W$, on note : $c(w) = BwB,\  ^wB = wBw^{-1}$.

$(A_3)$  $\forall s \in S, \forall w \in W$, $c(s).c(w)
\subset c(w) \cup c(sw)$.

$(A_4)$ $\forall s \in S$, $^sB \not\subset B$.}

A toute partie $S'$ de $S$, on fait correspondre le groupe $G_{S'} =
BW_{S'}B$ o\`u $W_{S'}$ est le sous-groupe de $W$ engendr\'e par
$S'$. $G_{S'}$ est le sous-groupe parabolique standard de type $S'$
et son rang est le cardinal de $S'$. Les sous-groupes paraboliques
maximaux standards sont donc  ceux de rang $\sharp S - 1$.

En particulier, soit $G$ un groupe de Lie semi-simple,
d'alg\`ebre de Lie $\got g$, $\got h$ une sous-alg\`ebre de Cartan
maximalement deploy\'ee de
$\got g$, $\got b$ une sous-alg\`ebre de Borel contenant $\got h$, $B$ un sous-groupe de
$G$ d'alg\`ebre de Lie $\got b$. Soit $\Delta$ un syst\`eme de racines
associ\'e et $\Pi$ une base de racines simples. 
Soit $M,M'$ respectivement
le centralisateur et le normalisateur  de $\got h$ dans $G$. Soit $S$
l'ensemble des reflexions associ\'ees aux racines simples. $S$
s'identifie \`a un sous-ensemble g\'en\'erateur du groupe de Weyl $W =
M'/M$, form\'e d'\'el\'ements d'ordre $2$. 
Dans ce cas, il est bien connu que le quadruplet $(G,B,M',S)$ est
un  syst\`eme de Tits.

Soit $G$ un groupe, $(G_i)$ une famille de
sous-groupes de $G$. 

\smallskip

{\bf D\'efinition 2.2 :} {\it On dit que $G$ est produit
amalgam\'e des $(G_i)$ suivant leurs intersections deux \`a deux si $G$ satisfait
\`a la propri\`et\'e universelle suivante :

Soit $H$ un groupe, $h_i : G_i \lra H$ une famille de morphismes
de groupes telle que :

$\forall x \in G_i \cap G_j, h_i(x) = h_j(x)$. 
Alors, il existe un et un seul morphisme de groupes $h:G \lra H$ tel
que : 
$\forall i, \forall x \in G_i, h_i(x) = h(x)$.}

Le r\'esultat suivant est d\^u \`a J.Tits (\cite{SE},
chapitre 2, I.7, corollaire 3).

\smallskip

{\bf Th\'eor\`eme 2.1 :} {\it Soit $(G,B,M',S)$ un syst\`eme de
Tits. On suppose que l'ensemble $S$ est de cardinal $n \geq 3$. Alors,
$G$ est produit amalgam\'e de ses sous-groupes paraboliques maximaux
standards suivant leurs intersections deux \`a deux.}

\section{\bf Les orbites nilpotentes sph\'eriques de  $\bf
sl_n({\mathbb R}), n \geq 4$.}

\subsection{\bf Quelques notations.} On adoptera, pour la suite, les notations suivantes : 

$\bullet$ Soit ${\got g} = sl_n({\mathbb R}), {\got g}_{\mathbb C} = sl_n({\mathbb
C})$ et $G$ un groupe de Lie connexe et simplement connexe d'alg\`ebre de 
Lie $\got g$. En fait, $G$ est le rev\^etement \`a deux feuillets de
$SL_n({\mathbb R})$. Soit $\mathcal K$ la forme de Killing d\'efinie sur
$\got g$, $\got
h$ une sous-alg\`ebre de Cartan d\'eploy\'ee de $\got g$, ${\got h}_{\mathbb C}$ la
sous-alg\`ebre de Cartan correspondante dans ${\got g}_{\mathbb C}$. Soit
$\Delta = \Delta ({\got g},{\got h}) = \{ \var_i - \var_j, i\not=j, 1
\leq i,j \leq n\}$ le syst\`eme de racines usuel pour $sl_n$, 
$\Delta^+ = \Delta ({\got g},{\got h}) = \{ \var_i - \var_j, 1
\leq i < j \leq n\}$, le syst\`eme de racines positives. Posons, de
plus, $\alpha_k = \var_k - \var_{k+1}, 1 \leq k \leq n-1$. Soit
$\Pi = \{ \alpha_k, 1 \leq k \leq n-1 \}$ le syst\`eme de racines
simples choisi. A chaque racine $\alpha$ dans $\Delta^+$, 
on associe le syst\`eme de
Chevalley usuel $(X_\alpha,H_\alpha, X_{-\alpha}), 
W_\alpha = X_\alpha - X_{-\alpha}$ et on pose, pour tout r\'eel $t$ :
 $$x_\alpha(t) = \exp_G tX_\alpha, \ x_{-\alpha}(t) = \exp_G tX_{-\alpha}$$
$$w_\alpha(t) = \exp_G tW_\alpha, h_\alpha(t) = \exp_G \ln |t|H_\alpha \ 
(t \not= 0)$$
$$w_\alpha^2 = \exp_G\pi W_\alpha$$
Soit $\pi_G : G \lra SL_n({\mathbb R})$ la projection
canonique correspondante et soit $z$ l'\'el\'ement non trivial du
noyau de $\pi_G$. 
Pour toute racine $\alpha$, on a : $w_\alpha^4 = z$.
Enfin, on d\'esignera par $\Gamma_\alpha$ le sous-groupe fini de $G$ engendr\'e 
par l'\'el\'ement $w_\alpha^2$.

$\bullet$ A  tout
sous-ensemble $(\alpha_i, \alpha_{i+1}, \dots, \alpha_{i+j})$ de
$\Pi$, on associe la sous-alg\`ebre de ${\got g}$ isomorphe \`a
$sl_{j+2}({\mathbb R})$, ayant  $(\alpha_i, \alpha_{i+1}, \dots,
\alpha_{i+j})$ comme syst\`eme de 
racines simples, et on la notera  $sl_{j+2}(\alpha_i, \dots,
\alpha_{i+j})$. Consid\'erons, pour $\dis{\var = \pm, 1 \leq k \leq [{n
\over 2}]}$,  la famille de vecteurs
suivante :
$$X_{\alpha_1}+ X_{-{\alpha_{n-1}}}, \dots,
 X_{\alpha_{k-1}}+ \var X_{-{\alpha_{n-k+1}}}, X_{-{\alpha_1}} +
X_{\alpha_{n-1}}, \dots, X_{-{\alpha_{k-1}}}+ \var X_{\alpha_{n-k+1}}$$
$$H_{\alpha_1}- H_{\alpha_{n-1}}, \dots,H_{\alpha_{k-1}}-
H_{\alpha_{n-k+1}}$$
Cette famille engendre une sous-alg\`ebre, isomorphe \`a $sl_k({\mathbb
R})$,  que nous noterons:
 $$sl_k(X_{\alpha_1}+ X_{-{\alpha_{n-1}}}, \dots,
 X_{\alpha_{k-1}}+ \var X_{-{\alpha_{n-k+1}}})$$ 

$\bullet$ Plus g\'en\'eralement, on notera $<X_1,\dots,X_p>$ le
sous-espace vectoriel de $\got g$, engendr\'e par la famille de vecteurs $X_1,\dots,X_p$.

$\bullet$ On introduit les  sous-alg\`ebres suivantes :
$${\got n} = \bigoplus_{\alpha \in \Delta^+} {\mathbb R}X_\alpha, \ 
{\got n}^- = \bigoplus_{\alpha \in \Delta^+} {\mathbb R}X_{-\alpha}, \ 
{\got b} = {\got h} \oplus {\got n}$$
$\got b$ est la sous-alg\`ebre de Borel associ\'ee \`a ce choix de racines positives et $B$ le sous-groupe
de Borel correspondant dans $G$. Soit $A = \exp {\got h}, N = \exp
{\got n}, N^- = 
\exp {\got n}^-$. 

$\bullet$ Soit enfin $K$ un sous-groupe compact maximal de $G$ et
$M$ le centralisateur de $A$ dans $K$. On sait que $M$ 
est un sous-groupe fini, engendr\'e par les $w_{\alpha}^2, \alpha \in 
\Delta^+$. On sait, \'egalement, que $B = M.A.N$. 

\subsection{}Les orbites
nilpotentes de ${\got g}_{\mathbb C}$ sont identifi\'ees aux partitions
de l'entier $n$. Selon la classification de D.Panyushev donn\'ee dans
\cite{PA}, les orbites
sph\'eriques correspondent aux partitions de la forme suivante :
$$(2^k, 1^{n-2k}), 1 \leq k \leq \frac{n}{2}$$  
On notera $O_k^C = (2^k,1^{n-2k})$ l'orbite correspondante et on
l'appellera orbite sph\'erique ``d'ordre $k$'' de ${\got g}_{\mathbb C}$. 

L'orbite $O_1^C$ \'etant l'orbite minimale, on ne consid\'erera donc que les
orbites sph\'eriques $O_k^C, 2 \leq k \leq \frac{n}{2}$.

Pour tout entier $k, 2 \leq k <\frac{n}{2}, \ O_k^C \cap \got g$ est
constitu\'e d'une seule orbite r\'eelle. Par contre, si
$n = 2p$, 
$O_p^C \cap {\got g}$ est r\'eunion de deux orbites r\'eelles.

D'apr\`es \cite{PA}, chaque $G$-orbite sph\'erique a pour
g\'en\'erateur une somme de vecteurs radiciels associ\'es \`a des
racines  simples deux \`a deux
orthogonales. 
Pour tout entier $\dis{k, 2 \leq k \leq  \frac{n}{2}}$ et pour $\var =
\pm 1$ , on pose :
$$Y_{k,\var} = \sum_{i=0}^{i=k-2} X_{\alpha_{2i+1}} + \var
X_{\alpha_{2k-1}}$$

Soit $O_{k,\var} = G.Y_{k,\var}$ l'orbite nilpotente r\'eelle
correspondante.  On constate que :

- Si $\dis{k <\frac{n}{2}}$, $O_{k,1}
= O_{k,-1}$ et on obtient ainsi une unique orbite nilpotente sph\'erique
r\'eelle d'ordre $k$. 

- Si $n = 2p$, $O_{p,1}$ et $O_{p,-1}$
sont deux orbites distinctes qui constituent les deux orbites
nilpotentes sph\'eriques r\'eelles d'ordre $p$.

On note $I_n = \{(k,\var), k \in ({\mathbb N} \cap [2,\frac{n}{2}]), 
 \var \in \{-1,1\} \}$. L'ensemble  
 $\{O_{k,\var}, (k,\var) \in I_n\}$ est l'ensemble des 
orbites nilpotentes sph\'eriques r\'eelles non minimales de $\got g$.

Enfin, on peut calculer ais\'ement, selon \cite{CO}, corollaire
6.1.4, la dimension de chaque orbite et on obtient :
$$\forall (k,\var) \in I_n, \dim O_{k,\var} = 2k(n-k)$$
Les orbites d'ordre $[\frac{n}{2}]$ joueront un r\^ole
particulier dans la suite. Nous les appellerons ``orbites sph\'eriques
maximales''.

Les g\'en\'erateurs pr\'ecis\'es 
pr\'ec\'edemment ne sont pas cependant des g\'en\'erateurs d'une
$B$-orbite ouverte.

Soit $(k,\var) \in I_n$. Posons :
$$
\begin{array}{rl}
\forall i,j, 1 \leq i \leq j \leq n-1, \ \beta_{i,j} &= \sum_{s=i}^{s=j}
\alpha_s \\
 \ \forall i, 1 \leq i \leq k, 
 \beta_i &= \beta_{i,n-i} \\ 
\ X_{k, \var} &= \sum_{i=1}^{i=k-1} X_{-\beta_i}
+ \var X_{-\beta_k} 
\end{array}
$$
On a, pour tout $(k,\var)$ dans $I_n$, $O_{k,\var} = G.X_{k,\var}$. 
Soit ${\got
g}(X_{k,\var}) = {\got r}_{k,\var} \oplus \ ^u{\got g}(X_{k,\var})$ le
stabilisateur de $X_{k,\var}$ dans $\got g$,
o\`u ${\got r}_{k,\var}$
d\'esigne un facteur r\'eductif et $^u{\got g}(X_{k,\var})$ le radical
unipotent de ${\got g}(X_{k,\var})$. 

On adoptera, pour la suite, les notations suivantes : 
$$\begin{array}{rl}
\forall j, 1 \leq j \leq n-1, \ {\got l}_{j} &= 
sl_{n-2j}(\alpha_{j+1}, \dots, \alpha_{n-j-1}), \ {\rm si} \ j \leq  \frac{n}
{2} - 1 \\
&= 0, \ {\rm sinon} \\
\forall j, 2 \leq j \leq n-1, \ {\got v}_{j, \var} &= sl_{j}(X_{\alpha_1}+X_{-\alpha_{n-1}},
\dots, X_{\alpha_{j-1}}+\var X_{-\alpha_{n-j+1}}) \oplus <H_{\alpha_j}
- H_{\alpha_{n-j}}> \\
{\got u}(X_{k,\var}) &= < X_{-\beta_{i,j}}, \ 1 \leq i \leq
k \leq j  \leq n-1> \\
{\got v}(X_{k,\var}) &= < X_{-\beta_{i,j}}, \ k+1 \leq i \leq n-k \leq j \leq
n-1>, \   \ {\rm si} \ k < \frac{n}{2}\\
&= 0, \  \ {\rm sinon} 
\end{array}
$$
Le calcul nous donne :
$$\begin{array}{rl}
^u{\got g}(X_k) &= {\got u}(X_{k,\var}) \oplus {\got
v}(X_{k,\var})  \\
 {\got r}_{k,\var} &= {\got l}_k \oplus {\got v}_{k,\var}
\end{array}
$$
On d\'eduit de ceci que, pour tout $(k,\var)$ dans $I_n$,  
${\got b}(X_{k,\var}) = {\got b} \cap {\got r}(X_{k,\var})$
et ainsi :
$$\dim {\got b}(X_{k,\var}) = \frac{(n-2k-1)(n-2k+2)}{2} +k $$
D'o\`u :
 $\dim B.X_{k,\var} = 2k(n-k), \forall (k,\var) \in I_n$.
Ainsi, 
$B.X_{k,\var}$ est une $B$-orbite ouverte dans $O_{k,\var}$. 

\smallskip

{\bf Lemme 3.1 :} {\it Soit $(k,\var) \in I_n$. Alors, la $B$-orbite $B.X_{k,\var}$ est l'unique $B$-orbite 
dense contenue dans $O_{k,\var}$.}

\smallskip

{\bf Preuve :} Il suffit, pour d\'emontrer ce r\'esultat, de calculer le 
nombre de $B$-orbites ouvertes contenues dans l'espace sym\'etrique
$G / G(X_{k,\var})$. T.Matsuki dans\cite{MA}, proposition 1, \'etablit une formule 
permettant de d\'eterminer ce nombre. On obtient, en appliquant cette 
formule, le r\'esultat souhait\'e. On pourra se r\'ef\'erer, pour un
calcul analogue, 
\`a \cite{SA}, d\'emonstration de la proposition 2.2.

\qed

\subsection{}Soit $(k,\var) \in I_n, G(X_{k,\var})$ le stabilisateur dans $G$ de $X_{k,\var}$ 
et $R_{k,\var}$ un 
facteur r\'eductif de $G(X_{k,\var})$. 
On pose, pour toute la suite : $w^2 = w^2_{\beta_1}, \Gamma =
\Gamma_{\beta_1}$. On obtient :
$$R_{k,\var} = \Gamma.(R_{k,\var})_0$$
o\`u $(R_{k,\var})_0$ d\'esigne la composante neutre de $R_{k,\var}$.

Soit, d'autre part, $f_{k,\var} = {\mathcal K}(X_{k,\var}, .)$
l'\'el\'ement de
${\got g}^*$ associ\'e \`a $X_{k,\var}$. On pose, enfin :
$$Y(f_{k,\var}) = \{ \tau \in R_{k,\var}^{\got g}, \tau(e) = - Id\}$$
$$Adm_{k} = Adm_G(f_{k,\var})= \{ \tau \in Y(f_{k,\var}), d\tau = 0\}$$

\smallskip

{\bf Proposition 3.2 :} {\it 

1) Pour tout $(k,\var)$ dans $I_n$, l'orbite $O_{k,\var}$ est 
admissible. 

2) Le groupe $\Gamma$ s'identifie \`a un sous-groupe de l'extension  
$R_{k,\var}^{\got g}$, not\'e encore $\Gamma$, et l'on a :
$$R_{k,\var}^{\got g} = \Gamma.((R_{k,\var})_0)^{\got g}$$

3) L'extension $((R_{k,\var})_0)^{\got g}$ admet exactement deux 
composantes connexes et l'application $x \lra (x,1)$ identifie
$(R_{k,\var})_0$ \`a  
la composante neutre de $((R_{k,\var})_0)^{\got g}$.

4) L'ensemble $Adm_{k}$ se d\'ecrit de la mani\`ere
 suivante : Soit $\chi$ le caract\`ere de $((R_{k,\var})_0)^{\got g}$ d\'efini par :
$$\chi (x,t(x)) = 1, \  \forall x \in (R_{k,\var})_0, \chi(e) = -1$$  

- Si $n = 2p$, on a : $\forall j, t(w^2) = 1$ et on consid\`ere 
le sous-groupe des caract\`eres de $\Gamma$ qui prennent la valeur $1$ sur 
les \'el\'ements $(w^4,1)$. Ce sous-groupe s'identifie au 
groupe $\dis{\widehat{({\mathbb Z} / 2{\mathbb Z})}}$ des caract\`eres de 
$\dis{({\mathbb Z} / 2{\mathbb Z})}$ et l'on a :
$$Adm_{k} = \{ \rho \otimes \chi, \ \rho \in 
\widehat{({\mathbb Z} / 2{\mathbb Z})} \}$$
- Si $n = 2p+1$, on a : $\forall j, t(w^2) = i$ et on 
consid\`ere le sous-groupe des caract\`eres de $\Gamma$ qui prennent la 
valeur $1$ sur 
les \'el\'ements $(w^8,1)$. Ce sous-groupe s'identifie au 
groupe $\dis{\widehat{({\mathbb Z} / 4{\mathbb Z})}}$ des caract\`eres de 
$\dis{({\mathbb Z} / 4{\mathbb Z})^k}$. Soit  
$\dis{\widehat{({\mathbb Z} / 4{\mathbb Z})^-}}$ le sous-ensemble 
form\'e des caract\`eres qui prennent la valeur $-1$ sur les \'el\'ements 
$(w^4, -1)$. Alors, on a : 
$$Adm_{k} = \{ \rho \otimes \chi, \ \rho \in 
\widehat{({\mathbb Z} / 4{\mathbb Z})^-} \}$$}

{\bf Preuve :} La premi\`ere assertion est un r\'esultat de 
 J.Schwartz~\cite{SC}. Pour la suite, il nous faut pr\'eciser l'extension
m\'etaplectique $((R_{k,\var})_0)^{\got g}$ . 

Consid\'erons, pour cela, le sous-espace $L$ de $\got g$ engendr\'e par 
les vecteurs suivants :
$$X_{\var_i - \var_j}, 1\leq i \leq k, k+1 \leq  j \leq n$$
 On constate que $L$ s'identifie \`a un sous-espace Lagrangien orient\'e 
de ${\got g} / {\got g}(X_{k,\var})$. A tout $x$ de $R_{k,\var}$, on
associe donc  
l'orientation relative $e(L,x.L)$ des lagrangiens $L$ et $x.L$. On 
v\'erifie successivement les faits suivants :

- $\forall x \in (R_{k,\var})_0, e(L,x.L) = 1$.

- $e(L, w^2.L) = (-1)^{n-2k} =
(-1)^n$.

On d\'efinit donc, pour tout $x$ de $R_{k,\var}$, le nombre
complexe $t(x)$  par :
$$\begin{array}{rl}
t(w^2) &= 1, \ {\rm si} \ 
n=2p \\
&= i, \ {\rm si} \ n=2p+1 \\
\forall x \in (R_{k,\var})_0, t(x) &=1 
\end{array}
$$
On en d\'eduit que :

- Le groupe $\Gamma$ est isomorphe \`a un sous-groupe de
$R_{k,\var}^{\got g}$. L'image de l'\'el\'ement $w^2$, par cet isomorphisme,
est $(w^2, 1)$ si $n = 2p$, $(w^2, i)$ si $n=2p+1$.  
La deuxi\`eme assertion de la proposition s'en suit.

- On peut identifier $(R_{k,\var})_0$ \`a un sous-groupe de 
$((R_{k,\var})_0)^{\got g}$ par l'application $x \lra (x,1)$ et on obtient :
$$((R_{k,\var})_0)^{\got g} = (R_{k,\var})_0 \cup e.(R_{k,\var})_0$$
ce qui donne la troisi\`eme assertion.

Supposons $n=2p$. Comme, pour tout $j$, $(w^2,1)^2 \in (R_{k,\var})_0$, 
on ne doit consid\'erer que les caract\`eres de $\Gamma$ qui prennent la 
valeur $1$ sur ces \'el\'ements et la derni\`ere assertion s'en suit.

Si $n = 2p+1$, $(w^2,1)^4 \in (R_{k,\var})_0$. Dans ces 
conditions, on ne s'int\'eresse qu'aux caract\`eres de $\Gamma$ qui 
prennent la valeur $1$ sur ces \'el\'ements. Cet ensemble est bien un 
groupe isomorphe \`a  $\dis{\widehat{{\mathbb Z} / 4{\mathbb Z}}}$. 
D'autre part, par d\'efinition d'un param\`etre d'admissibilit\'e, on doit 
choisir les caract\`eres $\rho$ qui v\'erifient de surcroit la propri\`et\'e :
$$\rho(w^4,-1) = \rho(w^4,1).\rho(e)= -1$$  
Ceci nous donne le r\'esultat souhait\'e. On constate ainsi que
l'ensemble $Adm_{k}$ est, dans tous les cas, un ensemble \`a deux
\'el\'ements qui ne d\'epend pas du param\`etre $\var$, ce qui
justifie la notation choisie. 
\qed

\section{\bf Les restrictions aux paraboliques maximaux.} 

On sait, d'apr\`es le lemme 3.1, que si $P$ est un sous-groupe parabolique 
maximal de $G$, chaque orbite nilpotente sph\'erique $O_{k,\var}$ contient une et 
une seule $P$-orbite dense. L'objet de ce paragraphe est d'\'etudier les 
$P$-orbites obtenues et, notamment, de les caract\'eriser suivant la 
param\'etrisation de 1.2.

\subsection{}A chaque racine simple $\alpha_i$, on associe la 
sous-alg\`ebre parabolique maximale ${\got p}_i$, obtenue \`a partir du 
sous-ensemble de racines $\Pi \backslash \{\alpha_i\}$. On d\'efinit la 
d\'ecomposition de Langlands  ${\got p}_i = {\got m}_i \oplus {\got a}_i 
\oplus {\got n}_i$ de ${\got p}_i$, avec :
$$\begin{array}{rl}
{\got m}_i &= sl_i(\alpha_1,\dots, \alpha_{i-1}) \oplus sl_{n-i}, \ \ (
\alpha_{i+1},\dots,\alpha_{n-1}) \hfill  (sl_1(\dots) = 0)\\
{\got a}_i &= {\mathbb R}H_{\alpha_i} \\
{\got n}_i& = \ ^u{\got p}_i 
\end{array}
$$

\smallskip

{\bf Remarque 4.1:} On remarque le fait important suivant : Chaque radical 
${\got n}_i$ est ab\'elien.

\smallskip

Soit $A_i = \exp {\got a}_i, N_i = \exp {\got n}_i, M_i = \Gamma_{\alpha_i}.
(M_i)_0$ et $P_i = M_i.A_i.N_i$ le sous-groupe parabolique maximal de $G$ 
correspondant.

Soit $(k,\var) \in I_n, P_i(X_{k,\var}) = P_i \cap G(X_{k,\var})$ le stabilisateur dans $P_i$ de $X_{k,\var}$. 
On peut \'ecrire~: $P_i(X_{k,\var}) = R_{i,k,\var}.^uP_i(X_{k,\var})$,
o\`u $R_{i,k,\var} = M_i \cap 
R_{k,\var}$ est un facteur r\'eductif de $P_i(X_{k,\var})$. Soit ${\got r}_{i,k,\var}$ 
l'alg\`ebre de Lie de $R_{i,k,\var}$. 
Posons, pour tout $i, 1 \leq i \leq n-1$ : 
$$
\begin{array}{rl}
{\got l}_{i,k} &= {\got l}_k \cap {\got m}_i \\
{\got v}_{i,k, \var} &= {\got v}_{k,\var} \cap {\got m}_i 
\end{array}
$$
On obtient :
$${\got r}_{i,k,\var} = {\got l}_{i,k} \oplus {\got v}_{i,k,\var}$$
En particulier, on constate que, pour tout $i, \ n-k \leq i \leq n-1,
{\got r}_{i,k,\var} = {\got r}_{n-i,k,\var}$. 

Introduisons ensuite les groupes finis $\Gamma_{i,k}$, d\'efinis de la
mani\`ere suivante :
$$
\begin{array}{rl}
\forall i, k < i < n-k, \Gamma_{i,k} &= \Gamma.\Gamma_{\alpha_i} \\
\forall i, i \leq k \ {\rm ou} \ i \geq n-k, \Gamma_{i,k} &= \Gamma 
\end{array}
$$
On obtient, alors : $R_{i,k,\var} = \Gamma_{i,k}.(R_{i,k,\var})_0$. 

On note enfin  
$f_{i,k,\var}$ la restriction de $f_{k,\var}$ \`a la sous-alg\`ebre ${\got p}_i$. 
On identifie naturellement $\got g$ \`a son dual ${\got g}^*$, via la forme
de Killing, et on note $Res_i : {\got g} \lra {\got p}_i^*$ l'op\'eration de 
restriction, moyennant l'identification pr\'ec\'edente.

\smallskip

{\bf Proposition 4.1 :} {\it 

i) Pour tout $(k,\var)$ dans $I_n$, la $P_i$-orbite 
$P_i.X_{k,\var}$ est l'unique $P_i$-orbite ouverte dense contenue dans
$O_{k,\var}$. 

ii) L'application $Res_i$ induit un diff\'eomorphisme  de
$P_i.X_{k,\var}$ sur  $P_i.f_{i,k,\var}$.  

iii) Chaque $P_i$-orbite $P_i.f_{i,k,\var}$ est $P_i$-admissible.}

\smallskip

{\bf Preuve :}
La premi\`ere assertion est due au fait que chaque orbite $O_{k,\var}$ contient une 
et une seule $B$-orbite ouverte dense. 

Soit $b_{k,\var}$ la restriction de $f_{k,\var}$ \`a $\got b$. On v\'erifie que 
${\got b}(X_{k,\var}) = {\got b}(b_{k,\var})$ et on a l'inclusion \'evidente $B(X_{k,\var}) 
\subset B(b_{k,\var})$. Soit, maintenant, $b_{k,\var,1}$ la restriction de $b_{k,\var}$ 
\`a ${\got n}_1$. En fait, $b_{k,\var,1}$ est la restriction \`a ${\got n}_1$ de 
${\mathcal K}(X_{-\beta_1},.)$.
De la double inclusion $ B\cap R_{k,\var} \subset B(b_{k,\var}) \subset B(b_{k,\var,1})$ et du 
fait que $N_1(b_{k,\var,1}) = \{1\}$, on en d\'eduit que $B(b_{k,\var}) = B \cap R_{k,\var}$. D'o\`u
l'\'egalit\'e des stabilisateurs : $B(X_{k,\var}) = B(b_{k,\var})$. De ceci, on d\'eduit 
que l'application "restriction"  induit un diff\'eomorphisme de l'orbite 
$B.X_{k,\var}$ sur $B.b_{k,\var}$. Comme ces orbites sont, respectivement, ouvertes 
et denses dans $P_i.X_{k,\var}$ et $P_i.f_{i,k,\var}$, on obtient ii).

Il existe un isomorphisme naturel d'espaces symplectiques de ${\got g} /
{\got g}(X_{k,\var})$ sur ${\got p}_i / {\got p}_i(X_{k,\var})$. Cet isomorphisme permet 
d'identifier l'extension $R_{i,k,\var}^{{\got p}_i}$ \`a un sous-groupe de 
$R_{k,\var}^{\got g}$. Il s'en suit que la restriction \`a $R_{i,k,\var}^{{\got p}_i}$
 d'un \'el\'ement de $Adm_{k}$ est un param\`etre d'admissibilit\'e de 
l'or\-bi\-te $P_i.X_{k,\var}$, ce qui prouve iii).
\qed

\smallskip

Nous noterons, dor\'enavant, $Adm_{i,k}$ l'ensemble des
param\`etres d'admissibilit\'e de l'or\-bi\-te $P_i.X_{k,\var}$.

\subsection{}On suppose toujours $(k,\var)$ dans $I_n$. Soit $h_{i,k,\var}$ la 
restriction de $f_{i,k,\var}$ \`a ${\got n}_i$ et ${\got b}_{i,k,\var} = {\got p}_i(h_{i,k,\var})$. Introduisons les
alg\`ebres suivantes : 
$$
\begin{array}{rl}
\ {\got g}_{i,k} &= {\got l}_i, \ \forall i < k \\
&= {\got l}_{i,k}, \ \forall i, \ k \leq i \leq n-k \\
&= {\got l}_{n-i},  \ \forall i > n- k \\
{\got v}'_{i,k,\var} &= {\got v}_{i,1}, \ \forall i < k \\
&= {\got v}_{i,k,\var}, \ \forall i, \ k \leq i \leq n-k \\
&= {\got v}_{n-i,1}, \ \forall i, \ i > n-k \\
{\got r}'_{i,k,\var} &= {\got g}_{i,k} \oplus {\got v}'_{i,k,\var}
\end{array}
$$
Le calcul montre alors que, pour tout $i, 1 \leq i \leq n-1$ :

- ${\got b}_{i,k,\var} = {\got r}'_{i,k,\var} + \ ^u{\got p}_i(X_{k,\var}) 
+ {\got n}_i$.

- ${\got r}'_{i,k,\var}$ est un facteur r\'eductif de ${\got
b}_{i,k,\var}$.

On a, de plus : 
\begin{equation}
\forall i,j, \ \ 1 \leq i \leq j \leq k, \ \ {\got
g}_{j,k} \subset {\got g}_{i,k}
\end{equation}

\smallskip

{\bf Remarque 4.2 :} Soit $n_{i,k}$ le rang de la sous-alg\`ebre
r\'eductive ${\got g}_{i,k}$. On constate que~:

- Si $k < [\frac{n}{2}]$, alors, pour tout $i, i < k$, $n_{i,k} =
n-2i-1 \geq 3$.  

- Si l'on
consid\`ere les orbites sph\'eriques maximales $O_{[\frac{n}
{2}]}$, alors :

\qquad - $ n_{p,p} = 0, \ p = [\frac{n}{2}]$.

\qquad - Si $n = 2p, {\got g}_{p-1,p} = sl_2(\alpha_p)$ et
$n_{p-1,p} =1$.

\qquad - Si $n = 2p+1, {\got g}_{p-1,p} =
sl_3(\alpha_p,\alpha_{p+1})$ et $n_{p-1,p} = 2$.
 
\qquad - $\forall i \leq [\frac{n}{2}] - 2, n_{i,k} \geq 3$.

Ces  remarques joueront un r\^ole
important, par la suite, dans la construction des repr\'esentations.

On d\'efinit, ensuite, la 
forme lin\'eaire $g_{i,k,\var}$ sur ${\got p}_i$ par :
$$
\begin{array}{rl}
\forall x \in {\got m}_i \oplus {\got a}_i, \ g_{i,k,\var}(x) &= 0, 
\\
 \forall x \in {\got n}_i, g_{i,k,\var} (x) &= f_{i,k,\var}(x)
\end{array}
  $$

\smallskip

{\bf Proposition 4.2 :} {\it ${\got b}_{i,k,\var}$ est une sous-alg\`ebre de
type fortement unipotent relativement \`a $g_{i,k,\var}$. $g_{i,k,\var}$ est une forme 
de type unipotent.} 

\smallskip

{\bf Preuve :} On commence par montrer que ${\got b}_{i,k,\var}$ est 
co-isotrope relativement \`a $g_{i,k,\var}$. En effet, soit $Z \in {\got b}_{i,k,\var}
^{\perp_g}$ un \'el\'ement de l'orthogonal dans ${\got p}_i$ de 
${\got b}_{i,k,\var}$ relativement \`a  $g_{i,k,\var}$.
 On \'ecrit $Z$ sous la forme $Z = x + y, x \in {\got m}_i \oplus 
{\got a}_i, y \in {\got n}_i$. On a :
$$g_{i,k,\var}([u,Z]) = 0 = g_{i,k,\var}([u,x]) + g_{i,k,\var}([u,y]), \forall u \in 
{\got n}_i \subset {\got b}_{i,k,\var}$$
0r, ${\got n}_i$ est ab\'elien, donc: 
$$g_{i,k,\var}([u,y])= 0 = g_{i,k,\var}([u,x]) = f_{i,k,\var}([u,x]) = h_{i,k,\var}([u,x]) $$
Comme ceci est vrai pour tout $u \in {\got n}_i$, on en d\'eduit que $x$ 
appartient \`a ${\got b}_{i,k,\var}$. Il en est de m\^eme de $y$. D'o\`u :
$${\got b}_{i,k,\var}^{\perp_g} \subset {\got b}_{i,k,\var}$$
Ceci montre bien que ${\got b}_{i,k,\var}$ est coisotrope relativement \`a 
$g_{i,k,\var}$.

On a l'inclusion triviale ${\got p}_i(g_{i,k,\var}) \subset {\got b}_{i,k,\var}$ et 
on sait, \'egalement, que ${\got r}'_{i,k,\var}$ est un facteur r\'eductif
de ${\got b}_{i,k,\var}$ contenu dans $ {\got m}_i \oplus 
{\got a}_i$.
Soit $X \in {\got r}'_{i,k,\var}, Y = A + B \in {\got p}_i,
 A \in {\got m}_i \oplus {\got a}_i, B \in {\got n}_i$. On a :
$g_{i,k,\var}([X,Y]) = g_{i,k,\var}([X,B])$, par d\'efinition de
$g_{i,k,\var}$. 
Comme $X \in {\got p}_i(h_{i,k,\var}), g_{i,k,\var}([X,B]) = 0$. D'o\`u $X \in
{\got p}_i(g_{i,k,\var})$ et, donc, ${\got r}'_{i,k,\var} \subset {\got
p}_i(g_{i,k,\var})$. Ainsi, ${\got r}'_{i,k,\var}$ est un facteur
r\'eductif de ${\got p}_i(g_{i,k,\var})$, 
ce qui montre que l'alg\`ebre ${\got b}_{i,k,\var}$ est bien de type
fortement unipotent relativement \`a $g_{i,k,\var}$.

Il est imm\'ediat de constater, pour finir, que la forme $g_{i,k,\var}$
est de type unipotent, puique $g_{i,k,\var}$ s'annule sur ${\got
r}'_{i,k,\var}$.
\qed

\smallskip

Soit $\lambda_{i,k,\var}$ la restriction de $f_{i,k,\var}$ \`a ${\got
r}'_{i,k,\var}$. On constate que la forme $\lambda_{i,k,\var}$ est nulle sur
${\got v}'_{i,k,\var}$. On
identifiera donc $\lambda_{i,k,\var}$ \`a une forme, encore not\'ee $\lambda_{i,k,\var}$, sur
${\got g}_{i,k}$.  

On en d\'eduit finalement, suivant les
notations de 1.2., le corollaire suivant   :

\smallskip

{\bf Corollaire 4.3 :} {\it Pour tout $(k,\var)$ dans $I_n$, les param\` etres
 de Duflo pour la
$P_i$-orbite $P_i.X_{k,\var}$ sont les suivants :
$$
\begin{array}{rl}
\forall i, i < k, \ P_i.X_{k,\var} &= O_{g_{i,k,\var}, \lambda_{i,k,\var}} \\
\forall i, k \leq i \leq  n-k, \ P_i.X_{k,\var} &= O_{f_{i,k,\var},0} \\
\forall i, n-k < i, \ P_i.X_{k,\var} &= O_{g_{i,k,\var},
\lambda_{i,k,\var}}
\end{array}
$$
En particulier, la $P_i$-orbite $P_i.X_{k,\var}$ est de type unipotent si et
seulement si $k \leq i \leq n-k$.  }

\section{\bf Construction de repr\'esentations unipotentes sph\'eriques 
par la m\'ethode des orbites.}
 
On rappelle que le but de ce travail est de
construire par la m\'ethode de Duflo-Torasso des repr\'esentations
unitaires irr\'eductibles de $G$ associ\'ees aux orbites $O_{k,\var}$, dans
le sens rappel\'e au paragraphe 0., c'est- \`a- dire des
repr\'esentations $\pi$ telles que l'annulateur infinit\'esimal de
$\pi$ dans $U({\got g}_{\mathbb C})$ ait pour vari\'et\'e des z\'eros la
cl\^oture de Zariski de $O_{k,\var}$. On va tout d'abord construire
des  familles de repr\'esentations ``candidates'' et on montrera dans les
paragraphes suivants que ces familles r\'epondent bien au probl\`eme
pos\'e.

\subsection{}Soit $(k,\var)$ dans $I_n$. La premi\`ere \'etape est de construire les $P_i$-repr\'esentations
associ\'ees aux orbites $P_i.X_{k,\var}$, selon les param\'etrisations
de Duflo donn\'ees dans les paragraphes 1.1 et 1.2.

Soit ${\got u}_{i,k,\var}= \ ^u{\got p}_i(X_{k,\var}) + {\got n}_i$
 le radical unipotent de ${\got b}_{i,k,\var},\  U_{i,k,\var} = \exp {\got
u}_{i,k,\var}$. La restriction de
$g_{i,k,\var}$ \`a ${\got 
u}_{i,k,\var}$ est nulle sur $^u{\got p}_i(X_{k,\var})$ et \'egale \`a
$h_{i,k,\var}$
sur ${\got n}_i$, nous la noterons encore $h_{i,k,\var}$. La
repr\'esentation de Kirillov $T_{i,k,\var}$ associ\'ee \`a $h_{i,k,\var}$ est
dans ce cas le caract\`ere $\rho_{i,k,\var}$ de $U_{i,k,\var}$ d\'efini par~:
$$\rho_{i,k,\var}(\exp X) = e^{-2i\pi h_{i,k,\var}(X)}, \forall X \in {\got
u}_{i,k,\var}$$

Soit $R'_{i,k,\var}$ un facteur r\'eductif 
de $P_i(g_{i,k,\var})$, d'alg\`ebre de Lie ${\got r}'_{i,k,\var}$. On
peut d\'ecrire l'extension m\'etaplectique $(R'_{i,k,\var})^{{\got p}_i}$
de la mani\`ere suivante :  

Soit $(G_{i,k})_0$ le sous-groupe analytique de $G$ d'alg\`ebre de Lie
${\got g}_{i,k}$,
$V'_{i,k,\var}$ le sous-groupe analytique de $G$ d'alg\`ebre
de Lie ${\got v}'_{i,k,\var}$.  
Posons : $G_{i,k} = \Gamma_{i,k}.(G_{i,k})_0$. On a, dans ces conditions
:
\begin{equation}
(R'_{i,k,\var})^{{\got p}_i} = (G_{i,k})^{{\got p}_i}.V'_{i,k,\var}
\end{equation}

Soit, enfin : $B_{i,k,\var} = R'_{i,k,\var}.U_{i,k,\var}$. 

D'apr\`es 1.1, la famille de $P_i$-repr\'esentations cherch\'ee est
param\'etr\'ee par l'ensemble $Y(g_{i,k,\var})$. Compte-tenu de 1.2 et du
corollaire 4.3, il
est n\'ecessaire de relier chaque repr\'esentation choisie dans
$Y(g_{i,k,\var})$ au param\`etre $\lambda_{i,k,\var}$. Or, la forme
$\lambda_{i,k,\var}$ est nulle sur ${\got v}_{i,k,\var}$, pour toute
valeur de l'entier $i$ et, d'apr\`es le corollaire 4.3,  nulle partout
si $k \leq i \leq n-k$. On va donc consid\'erer les sous-ensembles
suivants de $Y(g_{i,k,\var})$ :
$$
\begin{array}{rl}
\forall i, i < k, Y_{i,k,\var} &= \{ \tau \otimes 1, \tau \in
(G_{i,k})^{{\got p}_i }, \tau(e) = -Id\} \\
\forall i, k \leq i \leq n-k, Y_{i,k,\var} &= Adm_{i,k} \\
\forall i > n-k, Y_{i,k,\var} &= Y_{n-i,k,\var} 
\end{array}
$$

On constate qu'un \'el\'ement quelconque de $Adm_{i,k}$ est
essentiellement caract\'eris\'e par sa restriction au sous-groupe
$\Gamma$. Si $\tau \in Adm_{i,k}, \tau' \in
Adm_{j,k}$, nous dirons, par abus de langage, que $\tau  =
\tau'$ si les restrictions de $\tau$ et $\tau'$ \`a
$\Gamma$ sont les m\^emes. 

Les extensions $R_{i,k,\var}^{'{\got p}_i}$ et $R_{i,k,\var}^{'{\got
u}_{i}}$
sont \'egales et
la repr\'esentation m\'etaplectique $S_{i,k,\var}$ est donn\'ee par la
formule simple suivante : $\forall (x,t(x)) \in R_{i,k,\var}^{'{\got p}_i},
\ S_{i,k,\var}(x,t(x)) = t(x)$. 

Soit $\tau_{i,k,\var} \otimes 1  \in Y(g_{i,k,\var})$. La formule (2)
s'applique et nous  
donne la repr\'esentation correspondante de $B_{i,k,\var}$ d\'efinie
par :
$\forall (x,t(x)) \in
G_{i,k}^{'{\got p}_i}, \forall s \in V'_{i,k,\var}, \  \forall y
\in U_{i,k,\var},$ 
\begin{equation} 
\tau_{i,k,\var} \otimes 1 \otimes  S_{i,k,\var}T_{i,k,\var}(xsy) =
t(x)\rho_{i,k,\var}(y)\tau_{i,k,\var}(x,t(x))
\end{equation}

On d\'efinit enfin l'ensemble suivant :
$$
\begin{array}{rl}
{\mathcal Y}_{k,\var} = \{ \tau_{k,\var} = (\tau_{i,k,\var}) \in
\prod_{i=1}^{i=n-1}Y_{i,k,\var} \ / \ & . \forall i, \ 1 \leq i < k,
\tau_{n-i,k,\var} = \tau_{i,k,\var} \cr
& . \forall i,j, \  k \leq i,j \leq n-k, \tau_{i,k,\var} =
\tau_{j,k,\var}\} 
\end{array}
 $$
Soit $\tau_{k,\var} = (\tau_{i,k,\var}) \in {\mathcal Y}_{k,\var}$. On d\'efinit la 
famille $(\pi^k_i(\tau,\var))_{1 \leq i \leq n-1}$ de
$P_i$-repr\'esentations de type Duflo, \`a l'aide de (3),  par :  
\begin{equation}
\forall i, \ 1 \leq i \leq n-1, \ \pi^k_i(\tau,\var)  =
\ind_{B_{i,k,\var}}^{P_i} \tau_{i,k,\var} \otimes 1 \otimes  S_{i,k,\var}T_{i,k,\var}
\end{equation}

\subsection{}Le probl\`eme maintenant est de savoir comment construire une
repr\'esentation unitaire irr\'eductible de $G$
 \`a partir de la donn\'ee d'une famille de $P_i$-repr\'esentations de type
Duflo d\'efinie sur les sous-groupes paraboliques maximaux de $G$. La
m\'ethode utilis\'ee est en fait valable pour  
un groupe semi-simple $L$ 
, de rang r\'eel $p \geq 3$, et utilise la notion de produit amalgam\'e
introduite dans le paragraphe 2. 

Soit $\got l$ l'alg\`ebre de Lie de $L$. 
On consid\`ere une sous-alg\`ebre de Cartan ${\got h}_l$ de $\got l$,
 une sous-alg\`ebre de Borel ${\got b}_l$ de $\got l$,
$B_L$ un sous-groupe de $L$ d'alg\`ebre de lie ${\got b}_l$. On
consid\`ere ensuite un sous-groupe compact maximal $K_L$ de $L$ et un
syst\`eme de racines simples associ\'e \`a ${\got b}_l$. Soit $M'_L$
le normalisateur de ${\got h}_l$ dans $K_L$, et $S_L$ l'ensemble des
reflexions des racines simples. On sait, d'apr\`es le paragraphe 2,
que le quadruplet $(L,B_L,M'_L,S_L)$ est un syst\`eme de Tits. 

On
peut  d\'efinir la famille $\dis{(Q_i)_{1 \leq i \leq p}}$ des sous-groupes 
paraboliques maximaux standards de $L$. Le th\'eor\`eme 2.1 et la
propri\`et\'e d'universalit\'e du produit amalgam\'e nous permettent 
d'\'enoncer le r\'esultat 
suivant, dont la d\'emonstration est identique \`a celle du th\'eor\`eme 
4.11 de \cite{TO2}.

\smallskip

{\bf Th\'eor\`eme 5.1 :} {\it Soit $L$ un groupe de Lie semi-simple 
de rang $p$ sup\'erieur ou \'egal \`a $3$, soit $(Q_i)$ la famille 
des sous-groupes paraboliques maximaux standards de $L$, associ\'ee au 
sous-groupe $B_L$. Soit $\pi_i, 1 \leq i \leq p$, une 
repr\'esentation unitaire irr\'eductible de $Q_i$. On suppose que :

(i) Pour tous $i,j, 1 \leq i \leq j \leq p$, les restrictions de
$\pi_i$ et $\pi_j$ \`a $Q_i \cap Q_j$ sont \'equivalentes.

(ii) pour tout $i$, la restriction de $\pi_i$ \`a $B_L$ est irr\'eductible.

Alors, il existe une repr\'esentation unitaire irr\'eductible et une 
seule $\pi$ de $L$ telle que : 
$$\forall i, \ \pi_{\mid Q_i} = \pi_i$$}
On notera dor\'enavant : $\pi = A(L,\pi_i)$ une telle
repr\'esentation.

\subsection{}
Notre travail consiste maintenant \`a consid\'erer une famille
$(\pi_i^k(\tau,\var))_{1 \leq i \leq n-1}$ donn\'ee par (7) et \`a
d\'eterminer les param\`etres $\tau_{k,\var} \in {\mathcal Y}_{k,\var}$
pour lesquels les deux conditions du th\'eor\`eme 5.1 sont satisfaites. Il
suffira ensuite ``d'amalgamer'' les $\pi_i^k(\tau,\var)$
pour obtenir les repr\'esentations ``candidates''de $G$.   

Soit $\tau_{k,\var} \in {\mathcal Y}_{k,\var}, \ \  i,j$ deux indices tels
que : $1 \leq i \not= j \leq n-1$. On adopte les notations suivantes :

${\got p}_{i,j} = {\got p}_i \cap {\got p}_j, {\got
n}_{i,j} = \ ^u{\got p}_{i,j}$.

$P_{i,j} = P_i \cap P_j, \  N_{i,j} = N_i \cap N_j$.

${\got b}_{i,j,k,\var} = {\got b}_{i,k,\var} \cap {\got p}_j =
{\got p}_{i,j}(g_{i,k,\var}) + {\got n}_i, \  
B_{i,j,k,\var} = B_{i,k,\var} \cap P_j$.

Soit, d'autre part,  $\pi^k_{i,j}(\tau,\var)$ la restriction de $\pi^k_i(\tau,\var)$ \`a 
$P_{i,j}$. 

On va tout d'abord d\'eterminer les
param\`etres $(\tau_{k,\var})$ pour lesquels 
$\pi^k_{i,j}(\tau,\var)$ et $\pi^k_{j,i}(\tau,\var)$ sont
\'equivalentes. Compte-tenu de la d\'efinition de l'ensemble
${\mathcal Y}_{k,\var}$, il
suffit de se restreindre au cas suivant :
$ i < j, \ 1 \leq i < k, \ i+1  \leq j \leq n-i-1$.
On se placera dor\'enavant dans cette situation.
 
$\bullet$ L'orbite $B.X_{k,\var}$ est dense dans $P_i.X_{k,\var}$, par hypoth\`ese. Donc,
$P_i(X_{k,\var}).B$ est un ouvert de Zariski de $P_i$. Comme $P_i(X_{k,\var})
\subset B_{i,k,\var}$, on en d\'eduit que $B_{i,k,\var}.P_{i,j}$ contient un
ouvert de $P_i$ dont le compl\'ementaire dans $P_i$ est de codimension
1. Il s'en suit que :
\begin{equation}
\pi^k_{i,j}(\tau,\var) = \ind_{B_{i,j,k,\var}}^{P_{i,j}} (\tau_{i,k,\var}
\otimes 1 \otimes
S_{i,k,\var}.T_{i,k,\var})_{\mid B_{i,j,k,\var}}
\end{equation}  
On a, pour les m\^emes raisons :
\begin{equation}
\pi^k_{j,i}(\tau,\var) = \ind_{B_{j,i,k,\var}}^{P_{i,j}} (\tau_{j,k,\var}
\otimes 1 \otimes
S_{j,k,\var}.T_{j,k,\var})_{\mid B_{j,i,k,\var}}
\end{equation}  

$\bullet$ L'
alg\`ebre ${\got p}_{i,j,k} = {\got g}_{i,k} \cap {\got p}_j$  est en fait la sous-alg\`ebre
parabolique maximale de ${\got g}_{i,k}$ associ\'ee \`a la racine
$\alpha_j$, de d\'ecomposition de L\'evi 
${\got p}_{i,j,k} =
{\got m}_{i,j,k} \oplus {\got n}_{i,j,k}$. 
 Soit $P_{i,j,k} = G_{i,k} \cap P_j$ de d\'ecomposition de Levi 
$P_{i,j,k} = M_{i,j,k}.N_{i,j,k}$, le sous-groupe de $G$
d'alg\`ebre de Lie
${\got p}_{i,j,k}$.

D'apr\`es (4), on a l'inclusion :
$G_{j,k} \subset G_{i,k}, \forall j , i+1 \leq j \leq n-i-1$ et on v\'erifie que :
$G_{j,k} \subset M_{ij,k}$.

$\bullet $ $\tau_{i,k,\var}$ est une repr\'esentation de $(G_{i,k})^{{\got
p}_i}$. Or, la proposition 3.2, qui caract\'erise une telle extension
m\'etaplectique, permet d'affirmer qu'en fait $\tau_{i,k,\var}$ est
enti\`erement d\'etermin\'e par sa restriction \`a $G_{i,k}$, identifi\'e \`a un
sous-groupe de $(G_{i,k})^{{\got p}_i}$. Nous noterons de la m\^eme
facon cette restriction.

Soit $\lambda_{i,j,k,\var}$ la restriction de $\lambda_{i,k,\var}$ \`a ${\got
p}_{i,j,k}$. Soit ${\got c}_{i,j,k} = {\got g}_{j,k} \oplus {\got
n}_{i,j,k}$ et $C_{i,j,k}= G_{j,k}N_{i,j,k}$ le sous-groupe correspondant de
$P_{i,j,k}$. On v\'erifie alors alors le fait important suivant, dont la
d\'emonstration essentiellement technique ne sera pas reproduite ici:

\smallskip

{\bf Lemme 5.2 :} {\it ${\got c}_{i,j,k}$ est une sous-alg\`ebre de
type fortement unipotent de ${\got p}_{i,j,k}$, relativement \`a
$\lambda_{i,j,k,\var}$.}

\smallskip

$\bullet $ On peut, donc, suivant le paragraphe 1.1, associer au couple
$(\tau_{j,k,\var}, \lambda_{i,j,k,\var})$ une $P_{i,j,k}$-repr\'esentation de type
Duflo selon la formule (3)., dite $P_{i,j,k}$-repr\'esentation de type
Duflo associ\'ee au couple $(\tau_{j,k,\var}, \lambda_{i,j,k,\var})$.

\smallskip

{\bf Proposition 5.3 :} {\it Soit $\tau_{k,\var} \in {\mathcal
Y}_{k,\var}$. On suppose que, pour tout $i < k$ et tout $j, i+1 \leq j \leq
n-i-1$, la restriction de $\tau_{i,k,\var}$ au sous-groupe parabolique
$P_{i,j,k}$ de $G_{i,k}$ est la $P_{i,j,k}$-repr\'esentation de type
Duflo associ\'ee au couple
$(\tau_{j,k,\var},\lambda_{i,j,k,\var})$. Alors, pour tous $i,j, \ 1
\leq i<j \leq n-1$, les
restrictions \`a $P_{i,j}$ des 
repr\'esentations $\pi_i^k(\tau,\var)$ et $\pi_j^k(\tau,\var)$ sont
\'equivalentes.}   

\smallskip
 
{\bf Preuve :}

$\bullet$ Introduisons l'alg\`ebre ${\got b}'_{i,j,k,\var} = {\got
p}_{i,j}(g_{i,k,\var}) + {\got n}_{i,j}, B'_{i,j,k,\var}$ le sous-groupe de
$P_{i,j}$ correspondant. Soit, enfin, $h_{i,j,k,\var}$ la restriction de
$g_{i,k,\var}$ \`a ${\got n}_{i,j}$. On d\'efinit de m\^eme ${\got
b}'_{j,i,k,\var} = {\got p}_{i,j}(g_{j,k,\var}) + {\got n}_{i,j}$,
$B'_{j,i,k,\var}$ et $h_{j,i,k,\var}$.

Soit $T_{i,j,k,\var}$ la repr\'esentation de Kirillov du
radical unipotent de $B_{i,j,k,\var}$, associ\'ee \`a $h_{i,j,k,\var}$,
$S_{i,j,k,\var}$ la repr\'esentation m\'etaplectique. On
note, de la m\^eme fa\c con, $S'_{i,j,k,\var}, T'_{i,j,k,\var}$ les objets
correspondants associ\'es au radical unipotent de
$B'_{i,j,k,\var}$. On d\'efinit aussi $T_{j,i,k,\var}, S_{j,i,k,\var}$
et on remarque que $S'_{j,i,k,\var} = S'_{i,j,k,\var}, T'_{j,i,k,\var}
= T'_{i,j,k,\var}$. 

$\bullet $ Il est, tout d'abord,  facile de v\'erifier que  
$^u{\got b}_{i,j,k,\var}$ (resp. 
 $^u{\got b}_{j,i,k,\var}$ ) est une
sous-alg\`ebre de  $^u{\got b}'_{i,j,k,\var}$ (resp. de 
 $^u{\got b}'_{j,i,k,\var}$),  co-isotrope relativement \`a
$h_{i,j,k,\var}$ (resp. $h_{j,i,k,\var}$).

En appliquant le proc\'ed\'e d'induction par \'etages \`a la
repr\'esentation induite $\pi^k_{j,i}(\tau,\var)$, on obtient~:
$$\pi^k_{j,i}(\tau,\var) = \ind_{B'_{j,i,k,\var}}^{P_{i,j}}
(\ind_{B_{j,i,k,\var}}^{B'_{j,i,k,\var}} (\tau_{j,k,\var} \otimes 1 \otimes
S_{j,k,\var}T_{j,k,\var})_{\mid
B_{j,i,k,\var}})$$
En tenant compte de ce qui pr\'ec\`ede, on peut appliquer 
la proposition 1.1. et on a, alors, l'\'equivalence suivante :
$$\ind_{B_{j,i,k,\var}}^{B'_{j,i,k,\var}} (\tau_{j,k,\var} \otimes 1
\otimes S_{j,k,\var}T_{j,k,\var})_{\mid
B_{j,i,k,\var}} \simeq \tau_{j,k,\var} \otimes 1 \otimes 
S'_{i,j,k,\var}.T'_{i,j,k,\var}$$ 
On applique encore une fois le proc\'ed\'e d'induction par
\'etages, ce qui nous donne :
\begin{equation}
\pi^k_{j,i}(\tau,\var) =
\ind_{B'_{i,j,k,\var}}^{P_{ij}} (\ind_{B'_{j,i,k,\var}}^{B'_{i,j,k,\var}}
(\tau_{j,k,\var} \otimes 1 
\otimes S'_{i,j,k,\var}.T'_{i,j,k,\var}))
\end{equation}

Utilisons maintenant la proposition 1.2 et l'hypoth\`ese selon
laquelle la restriction de $\tau_{i,k,\var}$ \`a $P_{i,j,k}$ est la
$P_{i,j,k}$-repr\'esentation de type Duflo associ\'ee au couple $(\tau_{j,k,\var},\lambda_{i,j,k,\var})$. 
On obtient, alors :
$$\pi^k_{j,i}(\tau,\var) = \ind_{B'_{i,j,k,\var}}^{P_{ij}} (\tau_{i,k,\var}
\otimes 1 \otimes 
S'_{i,j,k,\var}.T'_{i,j,k,\var})$$
On applique enfin une nouvelle fois la proposition 1.1 ce qui nous
donne :
$$
\begin{array}{rl} 
\pi^k_{j,i}(\tau,\var)
&= \ind_{B'_{i,j,k,\var}}^{P_{ij}} (\ind_{B_{i,j,k,\var}}^{B'_{i,j,k,\var}}
(\tau_{i,k,\var} \otimes 1 \otimes  
S_{i,k,\var}.T_{i,k,\var})_{\mid B_{i,j,k,\var}}) \\
&= \ind_{B_{i,j,k,\var}}^{P_{ij}} (\tau_{i,k,\var} \otimes 1 \otimes
S_{i,k,\var}.T_{i,k,\var})_{\mid B_{i,j,k,\var}} \\
&= \pi^k_{i,j}(\tau,\var)
\end{array}
$$
\qed 

{\bf Proposition 5.4 :} 
{\it Soit $\tau_{k,\var} \in {\mathcal
Y}_{k,\var}$. On suppose que, pour tout $i < k$ et tout $j, i+1 \leq j \leq
n-i-1$, la restriction de $\tau_{i,k,\var}$ au sous-groupe parabolique
$P_{i,j,k}$ de $G_{i,k}$ est la $P_{i,j,k}$-repr\'esentation de type
Duflo associ\'ee au couple $(\tau_{j,k,\var},\lambda_{i,j,k,\var})$. Alors, 
pour tout $i$, la restriction de 
$\pi^k_i(\tau,\var)$ \`a $B$ est irr\'eductible.} 

\smallskip

{\bf Preuve :} Compte-tenu de la proposition 5.3, il suffit de montrer
que la restriction de $\pi^k_k(\tau,\var)$ \`a $B$ est irr\'eductible. Or,
il est facile de constater que cette restriction est une
$B$-repr\'esentation de type Duflo, qui est donc bien irr\'eductible.
\qed  

\subsection{\bf Le cas des orbites sph\'eriques non maximales.}
On va s'int\'eresser, tout d'abord, aux orbites sph\'eriques non 
maximales et donc supposer que $ k < [{n \over 2}]$. Nous
conserverons les notations introduites auparavant mais pouvons
remarquer que, dans ce cas, le param\`etre $\var$ peut \^etre choisi
\'egal \`a $1$, nous ne le ferons donc plus apparaitre dans les
notations qui vont suivre. On se fixe, pour la suite, un
param\`etre d'admissibilit\'e $\chi_k \in Adm_k$.

Compte-tenu de (4), on dispose d'une suite d\'ecroissante de groupes
semi-simples : $G_{k,k} \subset G_{k-1,k} \subset \dots \subset
G_{1,k}$. L'id\'ee consiste donc \`a construire les param\`etres
$\tau_{i,k}, 1 \leq i \leq k$, par r\'ecurrence sur $i$.

$\bullet $ On consid\`ere tout d'abord le groupe $G_{k,k}$, dont l'alg\`ebre de
Lie est ${\got g}_{k,k} =
{\got l}_k = sl_{n-2k}(\alpha_{k+1},\dots,\alpha_{n-k-1})$. La forme
$\lambda_{k,k}$ est nulle, on peut donc poser :
$$\tau_{k,k} = \chi_k$$

$\bullet $ Consid\'erons ensuite le groupe $G_{k-1,k}$, d'alg\`ebre de
Lie donn\'ee par~: $${\got
g}_{k-1,k} =sl_{n-2k+2}(\alpha_k,\dots,\alpha_{n-k})$$ 
La forme $\lambda_{k-1,k,\var}$ est,
dans ce cas, celle qui correspond \`a l'orbite minimale de
${\got g}_{k-1,k}$. On va donc appliquer \`a cette situation la m\'ethode
de construction de Torasso.

Suivant la proposition 5.3, on veut d\'efinir le param\`etre
$\tau_{k-1,k}$ de telle sorte que sa restriction \`a chaque parabolique 
$P_{k-1,j,k}, \ k \leq j \leq n-k,$ soit la $P_{k-1,j,k}$-repr\'esentation
de type Duflo associ\'ee au couple $(\chi_k,\lambda_{k-1,j,k})$, soit :
  
$$
(\tau_{k-1,k})_{\mid P_{k-1,j,k}} = \pi^{\chi}_{k-1,j,k} =
\ind_{C_{k-1,j,k}}^{P_{k-1,j,k}} 
(\chi_k \otimes S^{\lambda}_{k-1,j,k}.T^{\lambda}_{k-1,j,k})
$$

o\`u $S^{\lambda}_{k-1,j,k}$ et $T^{\lambda}_{k-1,j,k}$ sont
respectivement la repr\'esentation m\'etaplectique et la repr\'e\-sen\-tation de
Kirillov du radical unipotent de $P_{k-1,j,k}$, associ\'ees \`a la forme
$\lambda_{k-1,j,k}$.

On est dans la situation classique d\'ej\`a d\'ecrite par P.Torasso dans \cite{TO2} et on v\'erifie, comme dans \cite{TO2}, que les hypoth\`eses du th\'eor\`eme
5.1 sont satisfaites par la  famille de repr\'esentations
($\pi^{\chi}_{k-1,j,k}, k \leq j \leq n-k$). Comme le
groupe $G_{k-1,k}$ est de rang plus grand que $3$, on peut appliquer le
th\'eor\`eme 5.1. et on pose :
$$
\tau_{k-1,k} = A(G_{k-1,k}, \pi^{\chi}_{k-1,j,k}, k \leq j
\leq n-k)
$$
Comme, d'apr\`es la proposition 3.2., l'extension m\'etaplectique
$(G_{k-1,k})^{{\got p}_{k-1}}$ est engendr\'ee par $G_{k-1,k}$ et $e$, on
peut donc d\'efinir la repr\'esentation $\tau_{k-1,k}$ de
$(G_{k-1,k})^{{\got p}_{k-1}}$ par :
$$
\tau_{k-1,k}(x,t(x)) = \tau_{k-1,k}(x), \forall x \in G_{k-1,k},
\tau_{k-1,k}(e) = -Id
$$

On peut ensuite proc\'eder par r\'ecurrence. Supposons construites les
repr\'esentations $\tau_{j,k}, i+1 \leq  j \leq k$,
et consid\'erons la repr\'esentation $\tau_{i,k}$ du groupe $G_{i,k}$.
En reprenant les notations et la  construction pr\'ec\'edente, on
d\'efinit, pour tout 
$j, \ i+1 \leq j \leq n-i-1$,  les
$P_{i,j,k}$-repr\'esentations de type Duflo suivantes :
\begin{equation} 
\pi^{\chi}_{i,j,k} = \ind_{C_{i,j,k}}^{P_{i,j,k}}
\tau_{j,k} \otimes S^{\lambda}_{i,j,k}.T^{l}_{i,j,k}
\end{equation}
 
On applique  encore une fois le th\'eor\`eme 5.1 et on pose  :
\begin{equation}
\tau_{i,k} = A(G_{i,k},\pi^{\chi}_{i,j,k}, i+1 \leq j \leq n-i-1)
\end{equation}
Finalement, la repr\'esentation $\tau_{i,k}$ de l'extension
m\'etaplectique $G_{i,k}^{{\got p}_i}$ est encore donn\'ee par :
$$
\tau_{i,k}(x,t(x)) = \tau_{i,k}(x), \forall x \in G_{i,k},
\tau_{i,k}(e) = -Id
$$

Finalement, on constate que, dans ce cas, le param\`etre  $\tau_{k} \in
{\mathcal Y}_k$ est enti\`erement d\'etermin\'e par le param\`etre
d'admissibilit\'e $\chi_k$.

\smallskip

{\bf Th\'eor\`eme 5.5 :} {\it Soit $O_{k}$ l'orbite nilpotente
sph\'erique d'ordre $k < [{n \over 2}]$ de $\got g$ et $\chi_k$ un
param\`etre d'admissibilit\'e de $O_{k}$. Soit $\tau_{k} \in
{\mathcal Y}_k$, d\'efinie par r\'ecurrence, \`a
partir de $\chi_k$, par les formules (11) et (12). Soit
$(\pi^k_i(\chi))$ la famille correspondante de $P_i$-
repr\'esentations de type Duflo, donn\'ee par la formule (7). 
Alors, il existe une et une seule repr\'esentation unitaire
irr\'eductible de $G$, $\pi_{\chi,k}$, telle que :
$$\forall i, (\pi_{\chi,k})_{\mid P_i} = \pi^k_i(\chi)$$}
On notera, pour $\dis{k < [\frac{n}{2}]}$ :
$${\mathcal R}_k = \{ \pi_{\chi,k}, \chi \in Adm_k \}$$

\subsection{}{\bf Le cas des orbites sph\'eriques maximales .}    
On consid\`ere dans ce paragraphe le cas des orbites d'ordre maximal
$O_{p,\var}, \ p = 
[{n \over 2}]$. On rappelle, \`a ce sujet, que si $n=2p$ est pair, il existe
deux orbites d'ordre $p, \ O_{p,\var}, \ \var = \pm 1$, si $n = 2p+1$, il existe 
une seule orbite d'ordre $p$.  On se donne une
famille de $P_i$-repr\'esentations de 
type Duflo $(\pi^p_i(\tau,\var)_{1 \leq i \leq n-1})$ d\'efinie par
la formule (7) et on va 
donc construire la suite $(\tau_{i,p,\var})$, comme dans 5.4., par 
r\'ecurrence sur $i$. La m\'ethode est la m\^eme, seule la premi\`ere \'etape est 
diff\'erente. En effet, on rappelle  que :
$$
\begin{array}{rl}
{\got g}_{p-1,p} &= sl_2(\alpha_p), \  {\rm si} \ n=2p, \\
{\got g}_{p-1,p} &= sl_3(\alpha_p, \alpha_{p+1}), \ {\rm si} \ n=2p+1 \\
\end{array}
$$

\smallskip

{\bf Le cas $\bf n=2p$.} 

Notons $\chi_{p,\var'} , \var' = \pm 1$, le param\`etre
d'admissibilit\'e d\'efini par :
$$\chi_{p,\var'}(w^2) = (-1)^{\frac{1-\var'}{2}}$$

$\bullet$ Le param\`etre $\tau_{p,p,\var}$ est tout d'abord choisi dans
l'ensemble $ \{\chi_{p,\var'}, \var' = \pm 1\}$.

$\bullet$  On sait que : 
$$G_{p-1,p} = \Gamma.\widetilde{SL_2(\alpha_p)}$$
o\`u
$\widetilde{SL_2(\alpha_p)}$ est le rev\^etement \`a deux feuillets de
$SL_2(\alpha_p)$. 

On va d\'eterminer le param\`etre $\tau_{p-1,p,\var}$ en prenant soin qu'il
satisfasse aux hypoth\`eses impos\'ees dans la proposition 5.3. Pour cela,
on consid\`ere le sous-groupe de Borel   
$B_{\alpha_p}$ de
$\widetilde{SL_2({\alpha_p})}$. 
Comme on a : $w^2_{\alpha_p} = w^2.a$ o\`u $a$ est un \'el\'ement
appartenant \`a la composante neutre de $R_{p,\var}$, on peut
\'etendre la d\'efinition des param\`etres d'admissibilit\'e
$\chi_{p,\var'}$ au groupe $\Gamma_{\alpha_p}$.

Le param\`etre $\tau_{p-1,p,\var}$ doit donc satisfaire \`a la propri\'et\'e  suivante :
$$(\tau_{p-1,p,\var})_{\mid B_{\alpha_p}} =  \ind_{\Gamma_{\alpha_p}.\exp
{\mathbb R}X_{\alpha_p}}^{B_{\alpha_p}}  \chi_{p,\var'} \otimes t_{\alpha_p,\var }$$
o\`u $t_{\alpha_p,\var}$ est le caract\`ere d\'efini par :
$\dis{t_{\alpha_p,\var}(\exp tX_{\alpha_p}) = e^{-2i\pi\var t}}$. 

On peut r\'ealiser cette repr\'esentation induite dans $L^2({\mathbb
R}^{+*})$ (voir \cite{SA}, 3.7) et la propri\'et\'e pr\'ec\'edente  
se traduit alors par les  relations suivantes :
$$
\begin{array}{rl}
\forall f \in L^2({\mathbb R}^{+*}), \tau_{p-1,p,\var}(\exp
uX_{\alpha_p}).f(t) &= e^{-i\pi \var ut^2}.f(t) \\
\tau_{p-1,p,\var}(\exp
uH_{\alpha_p}).f(t) &= e^{u \over 2}f(te^u) \\
\tau_{p-1,p,\var}(w^2_{\alpha_p}).f(t) &=
\chi_{p,\var'}(w^2_{\alpha_p})f(t) 
\end{array}
$$
On constate, en \'etudiant la classification du dual unitaire de
$\widetilde{SL_2({\mathbb R})}$, que seules les s\'eries discr\`etes et
limites de s\'eries discr\`etes peuvent satisfaire aux relations
d\'efinies pr\'ecedemment. 

Les s\'eries discr\`etes (et limites de s\'eries discr\`etes) de 
$\widetilde{SL_2({\mathbb R})}$ sont param\'etr\'ees par les demi-entiers
${\mathbb Z}^* / 2$. Nous noterons $\dis{\tau_{\frac{n}{2}}, n \in 
{\mathbb Z}^{*},}$ ces  repr\'esentations.

On
sait  r\'ealiser la restriction au Borel d'une s\'erie
discr\`ete $\tau_{\var \frac{n}{2}}, \var = \pm 1, n \in {\mathbb N}^*,$ dans $L^2({\mathbb R}^{+*})$ (voir \cite{SA}, proposition
3.2., ou \cite{FA-KO}, th\'eor\`eme 13.1.1.). 
Les formules de la repr\'esentation sont les suivantes :
$$
\begin{array}{rl}
\forall f \in L^2({\mathbb R}^{+*}), \ \tau_{\var \frac{n}{2}}(\exp
uX_{\alpha_p}).f(t) &= e^{-i\var \pi ut^2}f(t) \\
\tau_{\var \frac{n}{2}}(\exp
uH_{\alpha_p}).f(t) &= e^{\frac{u}{2}}f(te^u) \\
\tau_{\var \frac{n}{2}}(w^2_{\alpha_p})
.f(t) &= (-1)^{\frac{n}{2}}(f(t) 
\end{array}
$$
 
On consid\`ere l' ensemble suivant :
\begin{equation}
{\mathcal D}_{\var,\var'}= \{ \chi_{p,\var'} \otimes  
\tau_{\var \frac{n}{2}}, \ n \equiv (1 - \var') \ mod(4) \}
\end{equation}
On v\'erifie ais\'ement que chaque \'el\'ement de ${\mathcal D}_{\var,\var'}$
est une repr\'esentation irr\'eductible de
$\Gamma.\widetilde{SL_2(\alpha_p)}= G_{p-1,p}$.
On constate, alors, que l'hypoth\`ese de
la proposition 5.3 est satisfaite d\`es que :
$\tau_{p,p,\var} = \chi_{p,\var'}$ et que   
$\tau_{p-1,p,\var}$ appartient \`a ${\mathcal D}_{\var,\var'}$.

Etant donn\'e un \'el\'ement $\delta \in {\mathcal
D}_{\var,\var'}$,
On pose, donc :
\begin{equation}
\tau_{p-1,p,\var} = \delta, \ \ 
\tau_{p,p,\var} = \chi_{p,\var'} 
\end{equation}

$\bullet$ Pour construire $\tau_{p-2,p,\var}$ on utilise le m\^eme
proc\'ed\'e  
que celui du paragraphe pr\'ec\'edent. On pose :
$$
\begin{array}{rl}
\forall j, p-1 \leq j \leq p+1, \ \pi^{\delta}_{p-2,j,p} &=
 \ind_{C_{p-1,j,p}}^{P_{p-2,j,p}}
\tau_{j,p,\var} \otimes
S^{\lambda}_{p-2,j,p,\var}.T^{\lambda}_{p-2,j,p,\var} \\
\tau_{p-2,p,\var} &= A(G_{p-2,p}, \pi^{\delta}_{p-2,j,p},
p-1 \leq j \leq p+1)
\end{array}
$$
$\bullet$  On construit ensuite les param\`etres $\tau_{i,p,\var}, 1 \leq i
\leq p-3$, 
 par r\'ecurrence selon des formules analogues \`a (11) et (12) soit~:
\begin{equation}
\pi^{\delta}_{i,j,p} = \ind_{C_{i,j,p}}^{P_{i,j,p}}
\tau_{j,p,\var} \otimes
S^{\lambda}_{i,j,p,\var}.T^{\lambda}_{i,j,p,\var} 
\end{equation}
\begin{equation}
\tau_{i,p,\var} = A(G_{i,p},\pi^{\delta}_{i,j,p}, i+1 \leq j \leq
n-i-1)
\end{equation}

\smallskip

{\bf Th\'eor\`eme 5.6 :} {\it Soit $O_{p,\var}$ une orbite nilpotente
sph\'erique maximale ``paire''  de $\got g$, $(n=2p)$. Soit ${\mathcal
D}_{\var,\var'}, \var' = \pm 1,$ l'ensemble de repr\'esentations 
de $\Gamma.\widetilde{SL_2(\alpha_p)}$ d\'efinie par (13) et soit 
$\delta  \in {\mathcal D}_{\var,\var'}$. Soit $\tau_{k,\var} \in
{\mathcal Y}_{k,\var}$, d\'efinie par r\'ecurrence, \`a partir de $\delta$,
par les formules 
(14),(15) et (16), et 
$(\pi^p_i(\delta))$ la famille correspondante de 
$P_i$-repr\'esentations de type Duflo donn\'ee par la formule (7).  
Alors, il existe une et une seule repr\'esentation unitaire
irr\'eductible de $G$, $\pi_{\delta,p}$,
telle que :
$$\forall i, (\pi_{\delta,p})_{\mid P_i} =
\pi^p_i(\delta)$$}
On notera :
$${\mathcal R}_{2p} = \{ \pi_{\delta,p} \ / \
\delta \in  
{\mathcal D}_{\var,\var'}, \var' = \pm 1\}$$

\smallskip

{\bf Le cas $\bf n = 2p+1$}

Comme dans 5.4., le param\`etre $\var$ peut \^etre choisi \'egal \`a $1$
et nous ne le ferons plus apparaitre dans les notations qui vont suivre.

$\bullet$
Dans ce cas, on peut \'ecrire :$G_{p-1,p} =
\Gamma.\widetilde{SL_3(\alpha_p,\alpha_{p+1})}$, o\`u
$\widetilde{SL_3(\alpha_p,\alpha_{p+1})}$ est le rev\^etement \`a deux
feuillets de $SL_3({\mathbb R})$. 
Dans \cite{TO1},  P.Torasso a d\'ecrit les 
repr\'e\-sen\-ta\-tions unitaires irr\'eductibles du groupe
$\widetilde{SL_3(\alpha_p,\alpha_{p+1})}$ associ\'ees, selon le sens
usuel, \`a 
l'orbite minimale de $sl_3(\alpha_p,\alpha_{p+1})$.Consid\'erons, pour
cela, les caract\`eres $t_z$ d\'efinis sur $\Gamma_{\alpha_p +
\alpha_{p+1}}$ par :
$$t_z(w^2_{\alpha_p + \alpha_{p+1}}) = z, \  z = -1,+1,-i,+i$$
Ces caract\`eres d\'efinissent  les param\`etres d'admissibilit\'e de
l'orbite minimale. 
On sait qu'il 
existe une et une seule repr\'esentation de $\widetilde{SL_3({\mathbb
R})}$, associ\'ee \`a $t_z$, 
que nous noterons $\rho_z$, lorsque $z = \pm 1$ (\cite{TO1}, proposition VI.
12) ou lorsque $z = -i$ (\cite{TO1}, th\'eor\`eme VII.1). Par contre, il
n'en existe pas lorsque $z = i$ (\cite{TO1}, th\'eor\`eme VII.1).

Par ailleurs, on rappelle que l'orbite $O_{p}$ poss\`ede deux
param\`etres d'admissibilit\'e $\chi_{p,\var'}, \var' = \pm 1,$ que
l'on peut d\'efinir par les donn\'ees suivantes :
$$\chi_{p,1}(w^2_{\alpha_p + \alpha_{p+1}}) = i, \  
\chi_{p,-1}(w^2_{\alpha_p + \alpha_{p+1}}) = -i$$

$\bullet$ On v\'erifie que $\chi_{p,-1} \otimes
\rho_{-i}$ est bien une repr\'esentation irr\'eductible de
$\Gamma.\widetilde{SL_3(\alpha_p,\alpha_{p+1})}$. En suivant un raisonnement analogue \`a ce qui pr\'ec\`ede, il est
facile de voir que les conditions impos\'ees par la proposition 5.3
impliquent la situation suivante :
\begin{equation}
\tau_{p-1,p} = \chi_{p,-1} \otimes \rho_{-i}, \ \ 
\tau_{p,p} = \chi_{p,-1} 
\end{equation}
 
$\bullet$ Pour construire $\tau_{p-2,p}$, on utilise les
formules suivantes  :
$$
\begin{array}{rl}
\forall j, p-1 \leq j \leq p+2, \ \pi^{\rho}_{p-2,j,p} &= \ind_{C_{p-2,j,p}}^{P_{p-2,j,p}}
\tau_{j,p} \otimes S^{\lambda}_{p-2,j,p}.T^{\lambda}_{p-2,j,p} \\
\tau_{p-2,p} &= A(G_{p-2,p}, \pi^{\rho}_{p-2,j,p}, p-1
\leq j \leq p+2) 
\end{array} 
$$
La construction des param\`etres $\tau_{i,p}, 1 \leq i \leq p-2$ proc\`ede ensuite 
par r\'ecurrence comme dans le cas 
pr\'ec\'edent, selon des formules analogues \`a (14) et (15) :
\begin{equation}
\pi^{\rho}_{i,j,p} = \ind_{C_{i,j,p}}^{P_{i,j,p}}
\tau_{j,p} \otimes
S^{\lambda}_{i,j,p}.T^{\lambda}_{i,j,p} 
\end{equation}
\begin{equation}
\tau_{i,p} = A(G_{i,p},\pi^{\rho}_{i,j,p}, i+1 \leq j \leq
n-i-1)
\end{equation}

\smallskip

{\bf Th\'eor\`eme 5.7 :} {\it Soit $O_{p}$ l'orbite nilpotente
sph\'erique maximale ``impaire''  de $\got g$, $(n=2p+1)$. Soit 
$\rho_{-i}$ la repr\'esentation minimale de
$\widetilde{SL_3({\mathbb R})}$ associ\'ee au param\`etre
d'admissibilit\'e $t_{-i}$ de l'orbite minimale  de $sl_3({\mathbb R})$. 
Soit $\tau_{k} \in
{\mathcal Y}_{k}$, d\'efinie par r\'ecurrence, \`a partir de $\rho_{-i}$,
par les formules 
(17),(18) et (19), et 
$(\pi^p_i(\rho))$ la famille correspondante de 
$P_i$-repr\'esentations de type Duflo donn\'ee par la formule (7).  
Alors, il existe une et une seule repr\'esentation unitaire
irr\'eductible de $G$, $\pi_{\rho,p}$,
telle que :
$$\forall i, (\pi_{\rho,p})_{\mid P_i} =
\pi^p_i(\rho)$$}

On notera :
$${\mathcal R}_{2p+1} = \{ \pi_{\rho,p}\}$$

Finalement, Les th\'eor\`emes 5.5, 5.6, 5.7 nous fournissent  des familles de
repr\'esentations unipotentes dont nous allons montrer maintenant
qu'elles sont bien associ\'ees aux orbites nilpotentes
sph\'eriques correspondantes. On pose :
$${\mathcal R} = \bigcup_{2 \leq k \leq \frac{n}{2}} {\mathcal R}_k $$

\section{\bf Sur la dimension de Gelfand-Kirillov des
\'el\'ements  de $\bf \mathcal R$.}

Soit $\pi$ une repr\'esentation dans $\mathcal
R$. On consid\`ere la repr\'esentation infinit\'esimale $\pi^\infty$
de $\pi$, qui est une repr\'esentation de l'alg\`ebre enveloppante
$U({\got g})$ de ${\got g}_{\mathbb C}$. Soit $I_\pi$ l'annulateur de $\pi^\infty$ dans $U({\got
g})$.

\smallskip

{\bf D\'efinition 6.1:} {\it On appelle Dimension de
Gelfand-Kirillov de $\pi$ la dimension de Gelfand-Kirillov de
l'alg\`ebre-quotient $U({\got g}) / I_\pi$. On note
$GKdim(\pi)$ cette dimension.}

\smallskip

{\bf D\'efinition 6.2:} {\it Soit $\pi \in {\mathcal R}$. $\pi$ sera dite
$GK$-associ\'ee \`a l'orbite $O_{k,\var}$ si l'on a :
$$GKdim (\pi) = \dim O_{k,\var}$$}

Notre but, dans ce paragraphe, est de d\'emontrer que chaque
\'el\'ement de ${\mathcal R}_k$ est $GK$-associ\'e \`a l'orbite $O_{k,\var}$.
Pour tout $i, \ 1 \leq i \leq n-1$, on consid\`ere la restriction 
de $\pi$ au parabolique maximal $P_i$. On pose : $I_{i,\pi} = I_\pi \cap
U({\got p}_i)$ et on peut
d\'efinir, comme pr\'ecedemment, la dimension de Gelfand-Kirillov,
$GKdim(\pi_{\mid P_i})$, de la repr\'esentation $\pi_{\mid P_i}$ par
la formule :
$$GKdim(\pi_{\mid P_i}) = GKdim (U({\got p}_i) / I_{i,\pi})$$
La premi\`ere \'etape consiste \`a relier $GKdim(\pi)$ aux nombres
$GKdim(\pi_{\mid P_i})$. 

\subsection{} On se place, dans ce paragraphe, dans le contexte plus
g\'en\'eral d'une alg\`ebre de Lie simple complexe $\got g$, de rang $r$
plus grand ou \'egal \`a $2$. On adopte, pour $\got g$, les
notations du paragraphe 3, pour une sous-alg\`ebre de Cartan, un
syst\`eme de racines ou encore les vecteurs-racine associ\'es. On
introduit \'egalement les sous-alg\`ebres paraboliques maximales
standard, ${\got p}_i, 1 \leq i \leq r$.

On pose : ${\got p}_{ij} = {\got p}_i \cap {\got p}_j, 1 \leq i
\leq j \leq r, \ {\got p}_{ii} = {\got p}_i$.

\smallskip

{\bf D\'efinition 6.3 :} 
{\it 

1) L'alg\`ebre de Lie $\got g$ est dite somme amalgam\'ee des
${\got p}_i$ suivant leurs intersections deux \`a deux si $\got g$ est
limite inductive du syst\`eme $({\got p}_{ij}, \varphi_{ij},
\varphi_{ji}, \ 1 \leq i
\leq j \leq r)$ o\`u $\varphi_{ij} : {\got p}_{ij} \lra {\got p}_i, \ 
\varphi_{ji} : {\got p}_{ij} \lra {\got p}_j$ sont les inclusions
canoniques. En d'autres termes, $\got g$ est solution du
probl\`eme universel suivant :

Etant donn\'ee une alg\`ebre de Lie $\got a$ et des morphismes
d'alg\`ebres de Lie $a_i : {\got p}_i \lra {\got a}$ tels que~:
$\forall (i,j), a_{i\mid {\got p}_{ij}} = a_{j\mid {\got p}_{ij}}$,
il existe un morphisme d'alg\`ebre de Lie et un seul $h : {\got g}
\lra {\got a}$ tel que : $h_{\mid {\got p}_i} = a_i, \ \forall i$.

2) L'alg\`ebre enveloppante  $U(\got g)$ est dite somme 
amalgam\'ee des
$U({\got p}_i)$ suivant leurs intersections deux \`a deux si $U(\got g)$ est
limite inductive du syst\`eme $(U({\got p}_{ij}), \phi_{ij},
\phi_{ji}, \ 1 \leq i
\leq j \leq r)$ o\`u $\phi_{ij} : U({\got p}_{ij}) \lra U({\got p}_i), \ 
\phi_{ji} : U({\got p}_{ij}) \lra U({\got p}_j)$ sont les inclusions
canoniques. En d'autres termes, $U(\got g)$ est solution du
probl\`eme universel suivant :

Etant donn\'ee une alg\`ebre associative  $\mathcal A$ et des morphismes
d'alg\`ebres  $u_i : U({\got p}_i) \lra {\mathcal A}$ tels que :
$\forall (i,j), u_{i\mid U({\got p}_{ij})} = u_{j\mid U({\got p}_{ij})}$,
il existe un morphisme d'alg\`ebres associatives 
  et un seul $h : U({\got g})
\lra {\mathcal A}$ tel que : $h_{\mid U({\got p}_i}) = u_i, \ \forall i$.}

\smallskip

{\bf Proposition 6.1 :} {\it Si $\got g$ est simple de rang au
moins \'egal \`a $2$, alors $\got g$ est somme amalgam\'ee de ses
paraboliques maximaux suivant leurs intersections deux \`a deux. De
m\^eme, $U({\got g})$ est somme amalgam\'ee des $U({\got p}_i)$
suivant leurs intersections deux \`a deux.}

\smallskip

{\bf Preuve :} 

1) Il suffit simplement de
prouver que $\got g$ satisfait \`a la propri\'et\'e universelle
donn\'ee dans la d\'efinition. Soit $\got a$ une alg\`ebre de lie et
$a_i : {\got p}_i \lra {\got a}$ des morphismes d'alg\`ebres de Lie
tels que :
$\forall (i,j), a_{i\mid {\got p}_{ij}} = a_{j\mid {\got p}_{ij}}$.
Soit ${\mathcal G} = \{ X_\alpha,H_\alpha,X_{-\alpha}, \alpha \in
\Pi\},$ o\`u $\Pi$ est un syst\`eme de racines simples dans $\got
g$. On sait que $\mathcal G$ est un ensemble de g\'en\'erateurs de $\got
g$. Or, puisque $\got g$ est suppos\'e de rang au moins $2$, pour
toute racine simple $\alpha$, il existe une sous-alg\`ebre parabolique
maximale ${\got p}_{i_\alpha}$ qui contient $X_{-\alpha}$. 
Posons :
$$\forall \alpha \in \Pi, \ u(X_\alpha) = a_{i_\alpha}(X_\alpha), 
u(H_\alpha) = a_{i_\alpha}(H_\alpha),
u(X_{-\alpha}) = a_{i_\alpha}(X_{-\alpha})$$
Compte-tenu des propri\'et\'es des morphismes $a_i$, 
on d\'efinit ainsi une application $u : {\mathcal G} \lra {\got a}$. Il
existe, alors, un morphisme d'alg\`ebres de Lie et un seul $h : {\got
g} \lra {\got a}$ tel que la restriction de $h$ \`a $\mathcal G$ soit
$u$. Ceci implique bien que $h_{\mid {\got p}_i} = a_i, \forall i$.

2) Posons : $U = U({\got g}), \forall i, U_i = U({\got p}_i),
\sigma : {\got g} \lra U, \ \sigma_i : {\got p}_i
\lra U_i$ les injections canoniques. 

Soit $\mathcal A$ une alg\`ebre associative et, pour tout indice
$i$, des morphismes d'alg\`ebres associatives  
$\beta_i : U_i \lra {\mathcal A}$ tels
que : $\forall (i,j), \ (\beta_i)_{ \mid U_i \cap U_j} = (\beta_j)_{ \mid
U_i \cap U_j}$. Soit
$\widetilde{\beta_i} = \beta_i \circ \sigma_i$. L'alg\`ebre associative
$\mathcal A$ \'etant  munie de sa structure usuelle d'alg\`ebre de Lie, il s'en
suit que, $\forall i,$ les morphismes $ \widetilde{\beta_i} : {\got
p}_i \lra {\mathcal A}$ sont des morphismes d'alg\`ebres de Lie satisfaisant \`a
la propri\'et\'e :
$$ \forall (i,j), (\widetilde{\beta_i})_{\mid {\got p}_{ij}} =
(\widetilde{\beta_j})_{\mid {\got p}_{ij}}$$
Compte-tenu de 1), il existe un morphisme d'alg\`ebres de Lie $h :
{\got g} \lra {\mathcal A}$ tel que : $\forall i, h_{\mid {\got p}_i} =
\widetilde{\beta_i}$. 

Par propri\'et\'e d'universalit\'e de l'alg\`ebre enveloppante
$U$, selon \cite{DI1}, lemme 2.1.3, on en d\'eduit un morphisme d'alg\`ebres
associatives $H : U \lra {\mathcal A}$ tel que :
$$H \circ \sigma = h, H(1) = 1$$
D\'esignons par $j_i : U_i \lra U$ l'injection canonique. Le fait que
$\forall i, H\circ \sigma_i = \widetilde{\beta_i}$ implique que :
$\forall i, H \circ j_i = \beta_i$, ce qui montre bien que $U$
satisfait \`a la propri\'et\'e universelle souhait\'ee.
\qed

\smallskip

En r\'ep\'etant les preuves pr\'ec\'edentes, on obtient le
corollaire suivant :

\smallskip

{\bf Corollaire 6.2 :} {\it Soit $\got g$ une alg\`ebre de Lie
semi-simple complexe de rang au moins \'egal \`a $2$. Soit ${\got p}_i, {\got
p}_j, i \not= j$ deux sous-alg\`ebres  paraboliques maximales de $\got
g$. Alors, $\got g$ (resp. $U({\got g})$) est somme amalgam\'ee de
${\got p}_i$ et ${\got p}_j$ (resp. $U({\got p}_i)$ et $U({\got
p}_j)$) suivant ${\got p}_{ij}$ (resp. $U({\got p}_{ij})$).}

\subsection{}Soit $\pi_k \in {\mathcal R}_k$ et, pour tout entier
$i, 1 \leq i \leq n-1, \ \pi_{i,k}$ la
restriction de $\pi_k$ au parabolique maximal $P_i$. Nous allons, tout d'abord, calculer la dimension de
Gelfand-Kirillov de $\pi_{k,k}$,
en rappelant que, d'apr\`es la construction faite dans le paragraphe 5,
cette repr\'esentation est la $P_k$-repr\'esentation de type Duflo associ\'ee au couple $(f_{k,k,\var}, \chi_{k})$, $f_{k,k,\var}$ \'etant
de type unipotent sur ${\got p}_k$ et $\chi_{k}$ \'etant un
param\`etre d'admissibilit\'e. 
 
Nous allons, pour cela, utiliser certains r\'esultats de \cite{DU} relatifs \`a la dimension de
Gelfand-Kirillov des repr\'esentations de type Duflo introduites dans
le paragraphe 1. En reprenant les notations de ce paragraphe, on
consid\`ere la  repr\'esentation $\pi_{q,\tau}$ d'un groupe presque
alg\'ebrique r\'eel $P$, o\`u $q$ est une forme de type unipotent sur
${\got p}$ et $\tau$ un \'el\'ement de $Y(q)$. On note $\mathcal E$ le
noyau de $\tau^\infty$ dans $U({\got p}(q))$.
Selon la classification des id\'eaux primitifs de M.Duflo, on fait
correspondre au couple $(iq,{\mathcal E})$ un id\'eal primitif de $U({\got
p})$, not\'e
$I_{iq,{\mathcal E}}$. On pourra se reporter \`a \cite{DU}, chapitre IV, pour la
d\'efinition d'un tel id\'eal. Les r\'esultats obtenus par M.Duflo
sont les suivants :

\smallskip

{\bf Proposition 6.3 :} {\it 

1) On a :
$$GKdim \ (U({\got p}) / I_{iq,{\mathcal E}}) = \dim ({\got p} /
{\got p}(q)) + GKdim \ (U({\got p}(q)) / {\mathcal E})$$

2) On a :
$$\ker \pi_{q,\tau}^\infty = \bigcap_{x \in P} x.I_{iq,{\mathcal E}}$$}
Consid\'erons, maintenant, le cas de la repr\'esentation 
$\pi_{k,k}$.
Par d\'efinition d'un param\`etre d'admissibilit\'e,
on constate que $U({\got p}_k(f_{k,k,\var})) / {\mathcal E}_k) = 0$, o\`u
${\mathcal E}_k$ est le noyau de $\chi_k$. 
D'autre part, soit $I_{\chi,k}$ l'id\'eal de $U({\got p}_k)$ correspondant au
couple $(if_{k,k,\var}, \chi_k)$. De la proposition 6.3, 2),  
on d\'eduit que~:
$$\ker \pi_{k,k}^\infty  = \bigcap_{x \in P_k} x.I_{\chi,k}$$
Comme $P_k$ ne poss\`ede qu'un nombre fini de composantes connexes,
l'intersection pr\'e\-c\'e\-den\-te ne se fait que sur un nombre fini de
termes et, suivant une propri\`et\'e classique des dimensions de
Gelfand-Kirillov, on a :
$$GKdim (U({\got p}_k) / \ker \pi_{k,k}^\infty) = \sup_{x \in P_k} 
(GKdim \ (U({\got
p}_k) / x.I_{\chi,k})) = GKdim \ (U({\got p}_k) / I_{\chi,k})$$
De la proposition 6.3, 1), et de ce qui pr\'ec\`ede, on d\'eduit le
lemme suivant :

\smallskip

{\bf Lemme 6.4 :} {\it On a :
$$GKdim \ \pi_{k,k} = \dim O_{k,\var}$$} 

Soit $I_{\pi,k}$ l'annulateur de
$\pi_k^\infty$ dans $U({\got g}), I_{\pi,i,k} = I_{\pi,k} \cap U({\got p}_{i,k})$. L'injection canonique $U({\got p}_k) \lra
U({\got g})$ induit un morphisme injectif de $U({\got p}_k) / I_{\pi,k,k}$
dans $U({\got g}) / I_{\pi,k}$ . On a donc l'in\'egalit\'e
:
\begin{equation}
GKdim(\pi_k) \geq \dim O_{k,\var} 
\end{equation}

On se propose, maintenant, de montrer le r\'esultat suivant :

\smallskip

{\bf Proposition 6.5 :} {\it 
Il existe une vari\'et\'e alg\'ebrique complexe ${\mathbb
X}_k$, dont la 
dimension est $\dis{\frac{1}{2}\dim O_{k,\var}}$, des morphismes d'alg\`ebres
$\phi_i : 
U({\got p}_i) \lra {\mathbb D}({\mathbb X}_k), i \in \{k,k+1\}$ ou $i \in \{k-1,k\}, {\mathbb
D}({\mathbb X}_k)$ d\'esignant 
l'alg\`ebre des op\'erateurs diff\'erentiels r\'eguliers sur ${\mathbb
X}_k$, tels que :
$$
\begin{array}{rl}
I_{\pi,i,k} &= \ker \phi_i, \\
(\phi_i )_{\mid \ U({\got p}_{i,j})} &=  
(\phi_j )_{\mid \ U({\got p}_{i,j})}, \ (i,j) \in \{k,k+1\} \
\hbox{ou} \ (i,j) \in \{k-1,k\} 
\end{array}$$}

On va d\'emontrer cette proposition en envisageant, tour \`a tour, le cas de 
l'orbite $O_{k}, 2 \leq k < [{n \over 2}]$ et le cas maximal de l'orbite 
$O_{[{n \over 2}],\var}$.

\smallskip

{\bf 1er Cas.} Supposons, tout d'abord,  que : $2 \leq k < [{n \over 2}]$. 
Dans ce cas, on sait que l'on peut faire disparaitre le param\`etre $\var$.

$\bullet $ On rappelle que :
$$\pi_{k,k} = \ind_{B_{k,k}}^{P_k} \chi_{k} \otimes 1 \otimes
S_{k,k}T_{k,k}$$
o\`u :
$${\got r}_{k,k} = {\got r}_{k}, \ \ {\got b}_{k,k} =
{\got r}_{k} + \ ^u {\got p}_k(X_{k}) + {\got n}_k$$
avec :
$$^u{\got p}_k(X_{k}) = < X_{-\beta_{ij}}, (i,j)  \in \{k+1,n-k\} \times  \{n-k,
n-1\}>$$

Soit
${\got q}_k = sl_{k+1}(\alpha_1,\dots,\alpha_k)$, $Q_k$ le sous-groupe
analytique de $G$ correspondant et ${\got q}_{k,k}$ la sous-alg\`ebre parabolique
maximale standard de ${\got q}_k$ associ\'ee \`a la racine
$\alpha_k$, $Q_{k,k}$ le sous-groupe parabolique de $Q_k$ correspondant. On note
$R_{Q,k}$ un facteur r\'eductif de $Q_{k,k}$, d'alg\`ebre de Lie
donn\'ee par : 
$$r_{Q,k} = sl_k(\alpha_1,\dots, \alpha_{k-1}) \oplus
{\mathbb R}H_{\alpha_k}$$

Soit, enfin, ${\got u}_k = \ ^u{\got p}_k(X_{k})^+ = < X_{\beta_{ij}},
X_{-\beta_{ij}} \in \ ^u{\got p}_k(X_{k})>$ et ${\got s}_k = r_{Q,k}
\oplus {\got u}_k$. 

${\got s}_k$ est une sous-alg\`ebre de ${\got
p}_k$ suppl\'ementaire de ${\got b}_{k,k}$ et de dimension $\dis{\frac{1}{2}\dim O_{k,\var}}$. On note ${\mathcal H}_{\chi,k}$ l'espace de Hilbert de la
repr\'esentation induite $\pi_{k,k}$ et on 
consid\`ere la vari\'et\'e r\'eelle ${\mathbb
X}_{k,{\mathbb R}} = R_{Q,k} \times {\got u}_k, {\mathbb X}_k$ sa
complexifi\'ee. Soit $j_k : {\mathbb X}_{k,{\mathbb R}}
\lra P_k$ l'application d\'efinie par : $j_k(x,Y) = x.\exp Y$. 
On d\'efinit ensuite l'espace de Hilbert
$E_k = L^2({\mathbb X}_{k,{\mathbb R}})$ 
des fonctions de carr\'e int\'egrable
sur ${\mathbb X}_{k, {\mathbb R}}$, muni de la mesure produit $dY \otimes
dx$, o\`u $dY$ est la mesure de Lebesgue sur ${\got u}_k$ et $dx$ la
mesure  de Haar sur $R_{Q,k}$.

Dans ces conditions, il est facile de voir que l'application $ {\got
j}_k :f \lra f \circ j_k$ d\'efinit une isom\'etrie de   
${\mathcal H}_{\chi,k}$ sur  $E_k$, donn\'ee par :
$$\forall x \in R_{Q,k}, \forall Y \in {\got u}_k, \forall f \in {\mathcal
H}_{\chi,k}, \ {\got j}_k(f)(x,Y) = f(x\exp Y)$$
A chaque $X \in {\got p}_k$, on associe l'op\'erateur diff\'erentiel
$l_X$, agissant sur l'espace $C^\infty (P_k)$, d\'efini par :
$$\forall f \in C^\infty(P_k), l_X.f(x) = {d \over dt}(f(\exp
-tXx))_{t=0}$$
L'application $X \lra l_X$ se prolonge en un morphismes d'alg\`ebres
de $U({\got p}_k)$ dans l'alg\`ebre des op\'erateurs diff\'erentiels
agissant sur $C^\infty(P_k)$. 

Notons ${\mathcal H}_{\chi,k}^\infty$ l'espace des vecteurs $C^\infty$
de $\pi_{k,k}$. Suivant le th\'eor\`eme 5.1 de \cite{PO} qui
caract\'erise les vecteurs $C^\infty$ d'une repr\'esentation induite,
on sait que :
$$
\begin{array}{rl}
 & - {\mathcal H}_{\chi,k}^\infty = \{ f \in C^\infty(P_k) \ / \ \forall a
\in U({\got p}_k), l_a.f \in {\mathcal H}_{\chi,k} \} \\
 & - \forall a \in U({\got p}_k), \forall f \in {\mathcal H}_{\chi,k}^\infty,
\ \pi_{k,k}^\infty(a).f = l_a.f 
\end{array}
$$
Par transport de structure, on r\'ealise $\pi_{k,k}$ dans
$E_k$ et on a :
$$E_k^\infty \subset C^\infty ({\mathbb X}_{k,{\mathbb R}})$$
 On peut ainsi
d\'efinir un morphisme d'alg\`ebres $\phi_k : U({\got p}_k) \lra {\mathbb
D}({\mathbb X}_k)$ tel que
:
\begin{equation}
\forall f \in E_k^\infty, \forall a \in U({\got p}_k),
\pi_{k,k}^\infty(a).f = \phi_k(a).f
\end{equation} 

$\bullet $ En ce qui concerne la repr\'esentation
$\pi_{k,k+1}$, on rappelle 
que :
$$\pi_{k,k+1}  = \ind_{B_{k+1,k}}^{P_{k+1}} \chi_{k} \otimes 1 \otimes
S_{k+1,k}T_{k+1,k,}$$ 
o\`u :
$${\got b}_{k+1,k} = {\got r}_{k+1,k} + \ ^u {\got p}_{k+1}(X_{k}) + {\got
n}_{k+1}$$ 
avec :
$${\got r}_{k+1,k} = sl_{n-2k-1}(\alpha_{k+2}, \dots, \alpha_{n-k-1})
\oplus {\got v}_{k,1} \oplus < H_{\alpha_{k+1}} >$$
$$
\begin{array}{rl}
^u{\got p}_{k+1}(X_{k}) = < X_{-\beta_{ik}}, i \in \{1,k\} > 
&\oplus <
X_{-\beta_{ij}}, (i,j) \in \{k+2, n-k\} \times \{n-k,n-1\} > \\
&\oplus < X_{\beta_{k+1j}}, k+1 \leq j \leq n-k-1 > 
\end{array}
$$
      
On pose : ${\got u}_{k+1} = (^u{\got p}_{k+1}(X_{k})\ \cap \ {\got
n}^-)^+ = < X_{\beta_{ij}}, \ X_{-\beta_{ij}} \in \ ^u{\got
p}_{k+1}(X_{k}) >$, puis ${\got s}_{k+1} = r_{Q,k} \ \oplus \ {\got
u}_{k+1}$.

${\got s}
_{k+1}$ est une sous-alg\`ebre de ${\got p}_{k+1}$ suppl\'ementaire de
${\got b}_{k+1,k}$ et de dimension $\dis{\frac{1}{2}\dim O_{k}}$. On
reprend ensuite les notations pr\'ec\'edentes :
${\mathcal H}_{\chi,k+1}$ est l'espace de Hilbert de la
repr\'esentation induite $\pi_{k,k+1}$,  
${\mathbb
X}_{k+1,{\mathbb R}} = R_{Q,k} \times {\got u}_{k+1}, {\mathbb X}_{k+1}$ est
sa complexifi\'ee. Soit $j_{k+1} : {\mathbb X}_{k+1, {\mathbb R}}
\lra P_{k+1}$ l'application d\'efinie par : $j_{k+1}(x,X) = x.\exp X$. 
Soit enfin $E_{k+1} = L^2({\mathbb X}_{k+1,{\mathbb R}})$
l'espace de Hilbert des fonctions de carr\'e
int\'egrable sur ${\mathbb X}_{k+1,{\mathbb R}}$, muni de la mesure
ad\'equate.
L'application $f \lra f \circ j_{k+1}$ est une isom\'etrie de ${\mathcal
H}_{\chi,k+1}$ sur $E_{k+1}$ et,
 selon les m\^emes
arguments que pr\'ecedemment, on peut d\'efinir un morphisme
d'alg\`ebres $\psi_{k+1} : U({\got p}_{k+1}) \lra {\mathbb D}({\mathbb
X}_{k+1})$ tel que :
\begin{equation}
\forall f \in E_{k+1}^\infty, \forall a \in U({\got p}_{k+1}),
(\pi_{k,k+1})^\infty(a).f = \psi_{k+1}(a).f 
\end{equation}

$ \bullet $ Les deux espaces ${\got u}_k$ et ${\got u}_{k+1}$
peuvent \^etre mis en dualit\'e par une transformation de Fourier
${\mathcal F}_k$, d\'efinie de la mani\`ere suivante :
A toute fonction $\varphi$ de l'espace de Schwartz ${\mathcal S}({\got
u}_{k+1})$, on associe la fonction ${\mathcal F}_k(\varphi)$, d\'efinie par
:
\begin{equation}
\forall Y \in {\got u}_k, {\mathcal F}_k(\varphi)(Y) = \int_{{\got
u}_{k+1}} \varphi (X) e^{2i\pi f_{k}([X,Y])} dX 
\end{equation}
Plus pr\'ecis\'ement, on pose $n_k = k(n-2k)$ et on identifie
naturellement les espaces ${\got u}_k$ et ${\got u}_{k+1}$ \`a ${\mathbb
R}^{n_k}$. L'application de ${\got u}_{k+1} \times {\got u}_k$ dans 
${\mathbb R}$ d\'efinie par : $(X,Y) \lra f_k([X,Y])$ s'identifie alors 
\`a une
forme bilin\'eaire $Q_k$ sur ${\mathbb R}^{n_k}$ donn\'ee par :
\begin{equation}
\forall X =(x_{ij}) \in {\got u}_{k+1}, \forall Y =(y_{ij}) \in
{\got u}_k, Q_k(X,Y) = \sum_{i=1}^{i=k} x_{ik}. y_{k+1,n-i} 
\end{equation}
L'application ${\mathcal F}_k$ se prolonge en un morphisme d'alg\`ebres,
not\'e encore ${\mathcal F}_k$, de ${\mathbb D}({\mathbb X}_{k+1})$ dans ${\mathbb
D}({\mathbb X}_k)$ et on pose : 
$$\phi_{k+1} = {\mathcal F}_k \circ
\psi_{k+1}$$
 $\phi_{k+1}$ est un morphisme d'alg\`ebres de $U({\got
p}_{k+1})$ dans ${\mathbb D}({\mathbb X}_k)$ et il est clair, d'apr\`es
$(21)$ et $(22)$, que $\ker \phi_k = I_{\pi,k,k}, \ker \phi_{k+1} =
I_{\pi,k+1,k}$. Il reste \`a prouver que la restriction de ces deux
morphismes \`a $U({\got p}_{k+1,k})$ est la m\^eme.

$\bullet $ Compte-tenu des d\'efinitions, il suffit de montrer que :

$\forall  Z \in \{X_{\alpha_k},
H_{\alpha_j}, k+1 \leq j \leq n-1, 
X_{-\alpha_{n-i}} \  1 \leq i \leq k, X_{\beta_{k+1,j}}, \ n-k \leq j
\leq n-1 \}$
$$\phi_k(Z) = \phi_{k+1}(Z)$$
 l'assertion \'etant imm\'ediate ou s'en d\'eduisant par
composition  pour les autres g\'en\'e\-ra\-teurs.

On obtient facilement les r\'esultats suivants :

\begin{equation}
\forall i, 1 \leq i \leq k, {\mathcal F}_k(\frac{\partial} 
{\partial x_{ik}}) = 2i\pi y_{k+1,n-i}
\end{equation}

\begin{equation}
\forall i, 1 \leq i \leq k, 2i\pi {\mathcal F}_k(x_{ik}) = - \frac{\partial} 
{\partial y_{k+1,n-i}}
\end{equation}

\begin{equation}
\forall i, 1 \leq i \leq k, {\mathcal F}_k(x_{ik}.\frac{\partial} 
{\partial x_{ik}}) = -Id -  y_{k+1,n-i}\frac{\partial}{\partial
y_{k+1,n-i}}
\end{equation}

\begin{equation}
\forall i,j,  1 \leq i < j \leq k, {\mathcal F}_k(x_{jk}.\frac{\partial} 
{\partial x_{ik}}) = - y_{k+1,n-i}.\frac{\partial}{\partial y_{k+1,n-j}}
\end{equation}

D'autre part, soit $Z \in {r_{Q,k}}$. On note $d_Z$ l'op\'erateur
diff\'erentiel sur $E_k^\infty$ d\'efini par :
$$\forall f \in E_k^\infty, \forall Y \in {\got u}_k, \forall x \in
R_{Q,k}, \  d_Z.f(x,Y) = {d \over dt}(f(x.\exp tZ.Y))_{t=0}$$

\smallskip

{\bf Lemme 6.6 :} {\it 

1) Il existe des fonctions $b_{ik}, \ C^\infty$ sur $R_{Q,k}$, 
telles que :
\begin{equation}
\phi_k(X_{\alpha_k}) = \phi_{k+1}(X_{\alpha_k}) 
= 2i\pi
\sum_{i=1}^{i=k} b_{ik}(.) y_{k+1,n-i}
\end{equation}

2)

$\forall i , 1 \leq i \leq k-1,$ \ 
\begin{equation}
\begin{split}
\phi_k(X_{-\alpha_{n-i}}) &= 
\phi_{k+1}(X_{-\alpha_{n-i}}) \\ 
&\displaystyle = d_{X_{\alpha_i}} +
\sum_{s=k+1}^{s=n-k} y_{s,n-i}\frac{\partial}{\partial
y_{s,n-i-1}}
\end{split}
\end{equation}

3)
\begin{equation}
\begin{split}
 \phi_k(H_{\alpha_{k+1}}) &= 
\phi_{k+1}(H_{\alpha_{k+1}}) \\
 &= -\sum_{i=1}^{i=k} y_{k+1,n-i}\frac{\partial}{\partial y_{k+1,n-i}} 
+ \sum_{i=1}^{i=k} y_{k+2,n-i}\frac{\partial}{\partial y_{k+2,n-i}} 
\end{split}
\end{equation}

4) $\forall j, k+2 \leq j \leq n-1,$
\begin{equation}
\begin{split}
 \phi_k(H_{\alpha_{j}}) 
&= \phi_{k+1}(H_{\alpha_{j}}) \\
 &= -\sum_{i=k+1}^{i=j-1} y_{i,j} \frac{\partial}{\partial y_{i,j}} 
-2y_{j,j}\frac{\partial}{\partial y_{j,j}} - \sum_{i=j+1}^{i=n-1} y_{j+1,i}
\frac{\partial}{\partial y_{j+1,i}} - {n-1 \over 2} Id 
\end{split}
\end{equation}

5)
\begin{equation}
\begin{split}
 \phi_k(X_{-\alpha_{n-k}}) 
&= \phi_{k+1}(X_{-\alpha_{n-k}}) \\ 
&= \sum_{i=k+1}^{i=n-k}\sum_{j=n-k}^{j=n-1} y_{i,n-k}y_{n-k,j}\frac{\partial} 
{\partial y_{i,j}} + \frac{n-k}{2}y_{n-k,n-k}   
\end{split}
\end{equation}

6)
Soit $a_{k+1,j}$ la fonction $C^\infty$ sur $R_{Q,k}$ 
d\'efinie par la relation :
$$x^{-1}x_{\beta_{k+1,j}}(u) x = x_{\beta_{k+1,j}}(u\ a_{k+1,j}(x))$$
Alors,
$\forall j , n-k \leq j \leq n-1,$
\begin{equation}
 \ \phi_k(X_{\beta_{k+1,j}}) = 
\phi_{k+1}(X_{\beta_{k+1,j}}) = -a_{k+1,j}(.)\frac{\partial}{\partial
y_{k+1,j}
} 
\end{equation}}

{\bf Preuve :} La d\'emonstration de ce lemme est essentiellement
technique. Nous nous contenterons d'en reproduire les calculs  pour
les formules (29) et (30).

1) Consid\'erons 
l'expression $x_{\alpha_k}(-u)x\exp Y, Y \in {\got u}_k$. Comme
$x_{\alpha_k}(-u)$ est 
un \'el\'ement du radical unipotent de $Q_{kk}$, sur lequel $R_{Q,k}$
agit, on a : $\dis{x_{\alpha_k}(-u)x = x. \prod_{i=1}^{i=k}
x_{\beta_{ik}}(c_{ik}(x,u))}$, o\`u $c_{ik} \in C^\infty(R_{Q,k} \times 
{\mathbb R})$ et $c_{ik}(x,0) = 0$. 

Ecrivons $\exp Y$ sous la forme :
$\dis{\exp Y = \prod_{k+1 \leq i
\leq n-k \leq j \leq n-1} x_{\beta_{ij}}(y_{ij})}$. 
On obtient : 
$$x_{\alpha_k}(-u)x\exp Y = x\exp Y \prod_{i=1}^{i=k} x_{\beta_{ik}}(c_{ik}
(x,u)).
\prod_{i=1}^{i=k}x_{\beta_{i,n-i}}(c_{ik}(x,u)y_{k+1,n-i})$$
Soit $f \in E_k$.

 Il s'en
suit que :
$$\pi_{k,k}(x_{\alpha_k}(u)).f = \prod_{i=1}^{i=k} e^{2i\pi
c_{ik}(x,u)y_{k+1,n-i}} .f$$
D'o\`u, par d\'erivation, et en posant :$\dis{\frac{\partial c_{ik}}
{\partial u}(x,0) = b_{ik}(x)}$, on a :
$$\phi_k(X_{\alpha_k}).f = 2i\pi \sum_{i=1}^{i=k} b_{ik}(.)y_{k+1,n-i}.f$$
On proc\`ede de la m\^eme fa\c con pour d\'eterminer
$\psi_{k+1}(X_{\alpha_k})$. Soit $X \in {\got u}_{k+1}$. On \'ecrit
$\exp X$ sous la forme : 
$$\exp X =
\prod_{i=1}^{i=k}x_{\beta_{ik}}(x_{ik}). \prod_{k+2 \leq i
\leq n-k \leq j 
\leq n-1} x_{\beta_{ij}}(x_{ij})$$
Ceci nous donne :
$$x_{\alpha_k}(-u)x\exp X = x \prod_{i=1}^{i=k}
x_{\beta_{ik}}(c_{ik}(x,u) + x_{ik}).\prod_{k+2 \leq i \leq n-k \leq j
\leq n-1} x_{\beta_{ij}}(x_{ij})$$
D'o\`u, par d\'erivation :
$$\psi_{k+1}(X_{\alpha_k}).f = \sum_{i=1}^{i=k} b_{ik}(.)
 \frac{\partial f}{\partial x_{ik}}$$
En utilisant la d\'efinition de $\phi_{k+1}$ et $(25)$, on en d\'eduit
$(29)$.

2) Soit $i, 1 \leq i \leq k-1$. Posons : $v_i(u) =
x_{\alpha_i}(u)x_{-\alpha_{n-i}}(u)$. On constate que $v_i(u) \in
B_{k,k,\var}$. En reprenant les notations
pr\'ec\'edentes, on peut \'ecrire, pour $Y \in {\got u}_k$ :

\begin{align*}
x_{-\alpha_{n-i}}(-u)x\exp Y  &= xx_{\alpha_i}(u)v_i(-u)\exp Y
\\
&=xx_{\alpha_i}(u)\prod_{s =
k+1}^{s=n-k}x_{\beta_{s,n-i-1}}(y_{s,n-i-1}+uy_{s,n-i}).\prod_{\substack{
k+1 \leq
s \leq n-k \\ n-k \leq j \leq n-1 \\ j \not= n-i-1}}x_{\beta_{sj}}(y_{sj}).
v_i(-u)
\end{align*}

D'o\`u, par d\'erivation :
$$\phi_k(X_{-\alpha_{n-i}}) =
d_{X_{\alpha_i}} +
\sum_{s=k+1}^{s=n-k} y_{s,n-i}\frac{\partial}{\partial
y_{s,n-i-1}}$$

En suivant le m\^eme type de calculs, on aboutit \`a la formule :
$$\psi_{k+1}(X_{-\alpha_{n-i}}) =
d_{X_{\alpha_i}} - x_{i+1,k}\frac{\partial}{\partial x_{ik}} +
\sum_{s=k+2}^{s=n-k} x_{s,n-i}\frac{\partial}{\partial
x_{s,n-i-1}}$$

On utilise ensuite la d\'efinition de $\phi_{k+1}$ et $(28)$ pour en
d\'eduire $(30)$.
A l'aide de calculs analogues et de $(25), (26), (27)$ et $(28)$,
on d\'emontre de m\^eme $(31), (32), (33)$ et $(34)$.

\smallskip

{\bf 2\`eme cas.} On suppose maintenant que $n = 2p$ et on envisage
le cas de l'orbite $O_{p,\var}$.

$\bullet $ Cette fois, la repr\'esentation $\pi_p$ d\'epend du choix
d'un param\`etre $\delta \in {\mathcal D}_{\var,\var'}$, ensemble  d\'efini par (13) . on va utiliser les sous-groupes
paraboliques $P_{p-1}$ et $P_p$. On rappelle que la forme
$f_p$ est de type unipotent, ce qui n'est pas le cas de $f_{p-1}$, et
l'on consid\`ere tout d'abord la repr\'esentation suivante :
$$\pi_{p,p} = \ind_{B_{p,p,\var}}^{P_p} \chi_{p,\var'} \otimes 1 \otimes
S_{p,p,\var}.T_{p,p,\var}$$
avec 
$${\got b}_{p,p,\var} = {\got r}_{p,\var} + {\got n}_p$$
On pose :
$$
\begin{array}{rl}
{\got u}_p &= < X_{\beta_{i,p-1}}, 1 \leq i\leq  p-1 > \oplus <
X_{\beta_{p+1,j}}, p+1 \leq j \leq 2p-1 > \\
{\got s}_p &= r_{Q,p-1}
\oplus {\mathbb R}H_{\alpha_p} \oplus {\got u}_p
\end{array}
$$

${\got s}_p$ est une sous-alg\`ebre de ${\got
p}_{p}$ suppl\'ementaire de ${\got b}_{p,p,\var}$ et de dimension
$\dis{\frac{1}{ 
2}\dim O_{p,\var}}$. On note ${\mathcal H}_{p}$ l'espace de Hilbert de la
repr\'esentation induite $\pi_{p,p}$ et on 
consid\`ere la vari\'et\'e r\'eelle ${\mathbb
X}_{p,{\mathbb R}} = R_{Q,p-1} \times {\mathbb R}^{+* }\times
{\got u}_p,  {\mathbb X}_p$ sa complexifi\'ee. Soit $j_p : {\mathbb
X}_{p,{\mathbb R}}
\lra P_p$ l'application d\'efinie par : 
$$\forall x \in R_{Q,p-1}, \forall t \in {\mathbb R}^{+*}, \forall Y \in
{\got u}_p, j_p(x,t,Y) = x.\exp Y\exp -\ln {t}H_{\alpha_p}$$ 
On consid\`ere ensuite l'espace de Hilbert $E_p = L^2({\mathbb
X}_{p,{\mathbb R}})$ et 
on d\'efinit l'application $ {\got
j}_p :   
{\mathcal H}_{p} \lra E_p$  par : 
$$\forall (x,t,Y) \in R_{Q,p-1} \times {\mathbb R}^{+*} \times {\got
u}_p, \forall f  \in {\mathcal H}_p, \ {\got j}_p(f)(x,t,Y)  = t^{{1
\over 2}-p} f \circ j_p(x,t,Y)$$
On constate qu'avec un bon choix de mesures  ${\got j}_p$ est 
une isom\'etrie.  
Par transport de structure, on r\'ealise $\pi_{p,p}$ dans
l'espace $E_p$ et
on a encore : $E_p^\infty \subset
C^\infty ({\mathbb X}_{p,{\mathbb R}})$. On utilise les arguments
d\'evelopp\'es dans le premier cas et le th\'eor\`eme 5.1. de \cite{PO}
pour d\'efinir un
morphisme d'alg\`ebres $\phi_p : U({\got p}_p) \lra {\mathbb
D}({\mathbb X}_p)$ tel que
:
\begin{equation}
\forall f \in E_p^\infty, \forall a \in U({\got p}_p),
\pi_{p,p}^\infty(a).f = \phi_p(a).f 
\end{equation}

$\bullet $ Consid\'erons maintenant la repr\'esentation :
$$\pi_{p-1,p} = \ind_{B_{p-1,p,\var}}^{P_p}
\delta \otimes 1 \otimes
S_{p-1,p,\var}.T_{p -1,p,\var}$$
avec 
$$
\begin{array}{rl}
{\got b}_{p-1,p,\var} &= {\got r}'_{p-1,p,\var} + \ ^u{\got
p}_{p-1}(X_{p,\var}) + {\got n}_{p-1} \\
{\got r}'_{p-1,p,\var} &= sl_2(\alpha_p) \oplus sl_{p-1}(X_{\alpha_1}+
X_{-\alpha_{2p-1}}, \dots, X_{\alpha_{p-2}}+X_{-\alpha_{p+2}}) \\
^u{\got p}_{p-1}(X_{p,\var}) &= <X_{-\beta_{pj}}, p \leq j \leq 2p-1> 
\end{array}
$$
On pose :
$$
\begin{array}{rl}
{\got u}_{p-1} &= <X_{\beta_{p,j}}, p+1 \leq j \leq 2p-1> \oplus 
<X_{\beta_{p+1,j}}, p+1 \leq j \leq 2p-1> \\
{\got s}_{p-1} &= r_{Q,p-1}
\oplus {\got u}_{p-1} 
\end{array}
$$

${\got s}_{p-1}$ est une sous-alg\`ebre de ${\got
p}_{p-1}$ suppl\'ementaire de ${\got b}_{p-1,p,\var}$. 
On note ${\mathcal H}_{\delta,p}$ l'espace de Hilbert de la
repr\'esentation induite $\pi_{p-1,p}$ et $V_\delta$ l'espace de
$\delta$.

On sait qu'il existe un op\'erateur unitaire $\phi_\delta : V_\delta
\lra L^2({\mathbb R}^{+*})$ permettant de r\'ealiser, dans l'espace $L^2({\mathbb
R}^{+*})$, la restriction de toute s\'erie
discr\`ete de $\widetilde{SL_2(\alpha_p)}$ \`a  son sous-groupe de Borel 
 (pour cela, on peut se r\'ef\'erer par exemple \`a
\cite{SA}, paragraphe 3.6). En particulier, on a les formules suivantes :
\begin{equation}
\begin{split}
\forall f \in L^2({\mathbb R}^{+*}), \delta(\exp uH_{\alpha_p}).f(t) &=
e^{\frac{u}{2}}.f(te^u) \\
\delta (\exp uX_{\alpha_p}).f(t) &=
e^{-i\pi\var ut^2}.f(t) 
\end{split}
\end{equation}
On 
consid\`ere la vari\'et\'e r\'eelle ${\mathbb
X}_{p-1,{\mathbb R}} = R_{Q,p-1} \times {\mathbb R}^{+*} \times
{\got u}_{p-1},  {\mathbb X}_{p-1}$ sa complexifi\'ee. Soit $E_{p-1} =
L^2({\mathbb X}_{p-1,{\mathbb R}})$. 

En utilisant l'op\'erateur $\phi_\delta$, on peut donc d\'efinir une
isom\'etrie ${\got j}_{p-1}$ de   
${\mathcal H}_{\delta,p}$ sur  $E_{p-1}$,
donn\'ee par :
$$\forall x \in R_{Q,p-1}, \forall X \in
{\got u}_{p-1}, \forall t \in {\mathbb R}^{+*}, \forall f \in {\mathcal
H}_{\delta,p}, \ {\got j}_{p-1}(f)(x,t,X) = \phi_\delta (f(x\exp X))(t)$$  
Par transport de structure, on r\'ealise $\pi_{p-1,p}$ dans
$E_{p-1}$ et on a encore : $E_{p-1}^\infty \subset
C^\infty ({\mathbb X}_{p-1,{\mathbb R}})$. On peut ainsi  
d\'efinir un morphisme d'alg\`ebres $\psi_{p-1} : U({\got p}_{p-1})
\lra {\mathbb D}({\mathbb X}_{p-1})$ tel que
:
\begin{equation}
\forall f \in E_{p-1}^\infty, \forall a \in U({\got p}_{p-1}),
\pi_{p-1,p}^\infty(a).f = \psi_{p-1}(a).f
\end{equation} 

$\bullet $ Comme pr\'ecedemment, on met en dualit\'e les espaces
${\got u}_p$ et ${\got u}_{p-1}$ par une transformation de Fourier
${\mathcal F}_p$ qui se d\'efinit selon une formule analogue \`a $(23)$,
soit :
\begin{equation}
\forall Y \in {\got u}_p, {\mathcal F}_p(\varphi)(Y) = \int_{{\got
u}_{p-1}} \varphi (X) e^{2i\pi f_{p,\var}([X,Y])} dX 
\end{equation}
On a \'egalement des formules analogues \`a (25), (26) (27) et (28), que 
l'on notera, sans les reproduire ici, (25'), (26'), (27') et (28'). 
Cette transformation  se prolonge en un morphisme d'alg\`ebres de ${\mathbb
D}({\mathbb X}_{p-1})$ dans ${\mathbb D}({\mathbb X}_p)$ et on pose enfin~:
$$\phi_{p-1} = {\mathcal F}_p \circ \psi_{p-1}$$
$\phi_{p-1}$ est un morphisme d'alg\`ebres de $U({\got
p}_{p-1})$ dans ${\mathbb D}({\mathbb X}_p)$ et il est clair, d'apr\`es
$(35)$ et $(36)$, que $\ker \phi_{p-1} = I_{\pi,p-1,p}, \ker \phi_p =
I_{\pi,p,p}$. Il reste \`a prouver que la restriction de ces deux
morphismes \`a $U({\got p}_{p-1,p})$ est la m\^eme.

Or, il est facile de constater, compte-tenu de la d\'efinition des
espaces de r\'ealisation consid\'er\'es, que les calculs donnant
$\phi_p (Z)$ ou $\phi_{p-1}(Z)$ sont les m\^emes que dans le premier
cas, en rempla\c cant $k$ par $p-1$, sauf pour $Z = H_{\alpha_p},
X_{\alpha_p}$.

On identifie l'espace ${\got u}_p$ \`a ${\mathbb
R}^{2p-2}$ \`a l'aide de l'application : $$Y \lra ((y_{i,p-1})_{1 \leq i \leq p-1},
(y_{p+1,j})_{p+1 \leq j \leq 2p-1}))$$
 On identifiera de la m\^eme
mani\`ere ${\got u}_{p-1}$ \`a ${\mathbb R}^{2p-2}$. On introduit, d'autre part, la
forme quadratique sur ${\mathbb R} \times {\mathbb R}^{2p-2}$, $\dis{{\mathcal
Q}_{p,\var}(t,Y) = \var t^2 - 
\sum_{i+j = 2p} y_{i,p-1}y_{p+1,j}}$. 

\smallskip

{\bf Lemme 6.7 :} {\it On a les
formules suivantes~:
\begin{equation}
\begin{split}
\phi_p(H_{\alpha_p}) &=
\phi_{p-1}(H_{\alpha_p}) \\
 &= (p-\frac{1}{2}).Id +
\sum_{i=1}^{i=p-1} y_{i,p-1} \frac{\partial}{\partial y_{i,p-1}} + 
\sum_{j=p+1}^{j=2p-1} y_{p+1,j} \frac{\partial}{\partial
y_{p+1,j}} + t\frac{\partial}{\partial t} 
\end{split}
\end{equation}
\begin{equation}
\phi_p(X_{\alpha_p}) =
\phi_{p-1}(X_{\alpha_p}) = 
i\pi{\mathcal Q}_{p,\var}.Id 
\end{equation}} 

{\bf Preuve :} Il s'agit encore de calculs essentiellement
techniques, on reproduira seulement ceux qui justifient $(39)$.

On
consid\`ere l'expression : $ h_{\alpha_p}(e^{-u}).x.\exp
Y.h_{\alpha_p}(t)$, pour $x \in R_{Q,p-1}, Y \in {\got u}_p, t \in
{\mathbb R}^{+*}$. On obtient :
$$h_{\alpha_p}(e^{-u}).x.\exp Y.h_{\alpha_p}(t) = 
x.\exp e^uY.h_{\alpha_p}(te^{-u})$$
En tenant compte de la d\'efinition de l'op\'erateur ${\got j}_p$ et
par d\'erivation, on aboutit \`a la formule souhait\'ee pour l'op\'erateur
$\phi_p(H_{\alpha_p})$.

Soit $X \in {\got u}_{p-1}$. On \'ecrit $\exp X$ sous la forme
suivante : 
$$\exp X = \prod_{j=p+1}^{j=2p-1} x_{\beta_{p,j}}(x_{p,j}).  
\prod_{j=p+1}^{j=2p-1} x_{\beta_{p+1,j}}(x_{p+1,j})$$ Un calcul identique au
pr\'ec\'edent nous donne : 
$$h_{\alpha_p}(e^{-u}).x.\exp X = x\prod_{j=p+1}^{j=2p-1} x_{\beta_{p,j}}(x_{p,j}e^{-u}).  
\prod_{j=p+1}^{j=2p-1} x_{\beta_{p+1,j}}(x_{p+1,j}e^u)h_{\alpha_p}(e^{-u})$$
Si l'on tient compte maintenant de la d\'efinition de l'op\'erateur
${\got j}_{p-1}$ et de son inverse, en utilisant $(36)$ et par
d\'erivation, on aboutit \`a :
$$\psi_{p-1}(H_{\alpha_p}) = 
 - \sum_{j=p+1}^{j=2p-1} x_{p,j} \frac{\partial}{\partial x_{p,j}} 
+ \sum_{j=p+1}^{j=2p-1} x_{p+1,j} \frac{\partial}{\partial
x_{p+1,j}} + t\frac{\partial}{\partial t} + \frac{1}{2}Id$$
On utilise ensuite la transform\'ee de Fourier ${\mathcal F}_p$ et la
formule $(27')$. On en d\'eduit le r\'esultat souhait\'e pour
$\phi_{p-1}(H_{\alpha_p})$.

\smallskip

{\bf 3\`eme cas.} On suppose enfin que $n=2p+1$ et on envisage le
cas de l'orbite maximale $O_{p}$. La repr\'esentation $\pi_p$ est
alors d\'etermin\'ee \`a partir du param\`etre $\rho_{-i}$ d\'efini dans le
paragraphe 5.5. 

$\bullet $ On consid\`ere tout d'abord la repr\'esentation : 
$$\pi_{p,p} = \ind_{B_{p,p}}^{P_p}
\chi_{p,-1} \otimes 1 \otimes S_{p,p}T_{p,p}$$
avec 
$${\got b}_{p,p} = {\got r}_{p} +\ ^u{\got p}_p(X_{p}) + {\got
n}_p$$
$$^u{\got p}_p(X_{p}) = < X_{-\beta_{p+1,j}}, p+1 \leq j \leq 2p
>$$
 
On note :
$$ 
\begin{array}{rl}
{\got u}_p &= <X_{\beta_{i,p-1}}, 1 \leq i \leq p-1> \oplus
<X_{\beta_{p+1,j}}, p+2 \leq j \leq 2p> \\ 
&\oplus <X_{\beta_{p+2,j}}, p+2
\leq j \leq 2p> \\ 
{\got d}_p &= <H_{\alpha_{p+1}},X_{\alpha_{p+1}}> \\
{\got s}_p &= r_{Q,p-1} \oplus {\got u}_p \oplus {\got v}_p 
\end{array}
$$

Alors, ${\got s}_p$ est une sous-alg\`ebre de ${\got p}_p$
suppl\'ementaire de ${\got b}_{p,p}$, de dimension $\dis{\frac{1}{2}
\dim O_p}$.

On note encore ${\mathcal H}_{p}$ l'espace de Hilbert de la
repr\'esentation induite $\pi_{p,p}$ et on 
consid\`ere la vari\'et\'e r\'eelle ${\mathbb
X}_{p,{\mathbb R}} = R_{Q,p-1} \times {\mathbb R}^* \times {\mathbb R}\times
{\got u}_p,  {\mathbb X}_p$ sa complexifi\'ee. Soit $j_p : {\mathbb X}_p
\lra P_p$ l'application d\'efinie par :
 
$\forall x \in R_{Q,p-1}, \forall (t,a) \in {\mathbb R}^{*} \times
{\mathbb R}, \forall Y \in
{\got u}_p,$
$$ j_p(x,t,a,Y) = x.\exp Y \exp aX_{\alpha_{p+1}}
w_{\alpha_{p+1}}^{1- {t \over |t|}}\exp \ln |t|H_{\alpha_{p+1}}$$ 

On d\'efinit ensuite, comme dans le cas pr\'ec\'edent, l'espace
 de Hilbert $E_p = L^2({\mathbb X}_{p,{\mathbb R}})$,
puis l'application $ {\got j}_p :   
{\mathcal H}_{p} \lra E_p$  par : ${\got j}_p = f \circ j_p$.

On constate que, pour un bon choix des mesures, ${\got j}_p$ est une 
isom\'etrie.  
Par transport de structure, on r\'ealise $\pi_{p,p}$ dans
l'espace $E_p$ et
on a encore : $E_p^\infty \subset
C^\infty ({\mathbb X}_{p,{\mathbb R}})$. On peut ainsi
d\'efinir un morphisme d'alg\`ebres $\phi_p : U({\got p}_p) \lra {\mathbb
D}({\mathbb X}_p)$ tel que
:
\begin{equation}
\forall f \in E_p^\infty, \forall a \in U({\got p}_p),
\pi_{p,p}^\infty(a).f = \phi_p(a).f 
\end{equation}

$\bullet $ Venons-en maintenant  \`a la repr\'esentation :  
$$\pi_{p-1,p}= \ind_{B_{p-1,p}}^{P_{p-1}}
\rho_{-i}\otimes 1 \otimes S_{p-1,p}T_{p-1,p}$$
avec 
$${\got b}_{p-1,p} = {\got r}'_{p-1,p} + \ ^u{\got p}_{p-1}(X_{p}) + {\got
n}_{p-1}$$
$${\got r}'_{p-1,p} = sl_3(\alpha_p,\alpha_{p+1}) \oplus
sl_{p-1}(X_{\alpha_1}+X_{-\alpha_{2p-1}}, \dots, X_{\alpha_{p-2}} +
X_{-\alpha_{p+2}}) \oplus <H_{\alpha_p} - H_{\alpha_{p+1}}>$$
$$
\begin{array}{rl} 
^u{\got p}_{p-1}(X_{p}) &= < X_{-\beta_{p+1,j}}, p+1 \leq j \leq 2p >
 \oplus <X_{-\beta_{p,j}}, p \leq j \leq 2p> \\ 
&\oplus <X_{\beta_{i,p-1}} +
X_{-\beta_{p+2,2p+1-i}}, 1 \leq i \leq p-1> 
\end{array}
$$

On rappelle \'egalement que le param\`etre $\rho_{-i}$ 
est celui
donn\'e par le th\'eor\`eme 5.7. et correspond \`a l'unique
repr\'esentation minimale de $sl_3$ associ\'ee au param\`etre
d'admissibilit\'e $t_{-i}$de l'orbite minimale. 
On note, ensuite : 
$$\begin{array}{rl}
{\got u}_{p-1} &= <X_{\beta_{p,j}}, p+2 \leq j \leq 2p> \oplus
<X_{\beta_{p+1,j}}, p+2 \leq j \leq 2p> \\
&\oplus <X_{\beta_{p+2,j}}, p+2
\leq j \leq 2p> \\  
{\got s}_{p-1} &= r_{Q,p-1} \oplus {\got u}_{p-1}
\end{array}
$$
Alors, ${\got s}_{p-1}$ est une sous-alg\`ebre de ${\got p}_{p-1}$
suppl\'ementaire de ${\got b}_{p-1,p}$.

On consid\`ere, dans un premier temps, la repr\'esentation minimale 
$\rho_{-i}$ de $\widetilde SL_3(\alpha_p,\alpha_{p+1})$.
On en donne un espace de r\'ealisation, selon des arguments d\'evelopp\'es 
par P.Torasso dans \cite{TO1}. Pour cela, on consid\`ere le sous-groupe 
parabolique $S_p$, associ\'e \`a la racine $\alpha_p$. On 
sait, alors, que la restriction de $\rho_{-i}$ \`a $S_p$ est une 
repr\'esentation induite d'un sous-groupe de $S_p$ dont l'alg\`ebre de Lie 
est une sous-alg\`ebre de type fortement unipotent relativement \`a 
l'orbite $S_p.X_{-\alpha_p - \alpha_{p+1}}$. On note $V_\rho$ l'espace 
de cette repr\'esentation induite. 

On d\'efinit, comme pr\'ecedemment, l'application ${\got j}_p : 
{\mathbb R}^* \times {\mathbb R} \lra S_p$ par :
$${\got j}_p(t,a) = \exp aX_{\alpha_{p+1}}.w_{\alpha_{p+1}}^{1 - 
\frac{t}{|t|}}.\exp \ln |t| H_{\alpha_{p+1}}$$
Cette application induit une isom\'etrie de 
$\phi_\rho : V_\rho \lra L^2({\mathbb R}^* \times {\mathbb R})$ qui permet
de r\'ealiser  
$\rho_{-i}$, par transport de structure, dans l'espace de Hilbert 
$L^2({\mathbb R}^* \times {\mathbb R})$.

Soit ${\mathcal H}_{\rho,p}$ l'espace de Hilbert de la
repr\'esentation induite $\pi_{p-1,p}$. On consid\`ere la
vari\'et\'e r\'eelle  ${\mathbb X}_{p-1, {\mathbb R}} = R_{Q,p-1} \times 
{\mathbb R}^* \times {\mathbb R} \times {\got u}_{p-1}$ et ${\mathbb X}_{p-1}$
sa complexifi\'ee. Soit $E_{p-1} = L^2({\mathbb X}_{p-1,{\mathbb R}})$. On
d\'efinit l'application ${\got j}_{p-1} : {\mathcal H}_{\rho,p} \lra
E_{p-1}$ par :
$$\forall f \in {\mathcal H}_{\rho,p}, \forall (x,t,a,X) \in {\mathbb
X}_{p-1,{\mathbb R}}, {\got j}_{p-1}(f)(x,t,a,X) = \phi_\rho(f(x\exp X))(t,a)$$
On a 
encore l'inclusion : $E_{p-1}^\infty \subset C^\infty({\mathbb X}_{p-1})$ et on 
peut ainsi d\'efinir un morphisme d'alg\`ebres $\psi_{p-1} : 
U({\got p}_{p-1}) \lra {\mathbb D}({\mathbb X}_{p-1})$ tel que :
\begin{equation}
\forall f \in E_{p-1}^\infty, \forall a \in U({\got p}_{p-1}),
\pi_{p-1,p}^\infty(a).f = \psi_{p-1}(a).f 
\end{equation}  

$\bullet $ Comme dans les cas pr\'ec\'edents, il existe une 
transformation de Fourier ${\mathcal F}_p$, d\'efinie par une formule 
analogue \`a (23) ou (38), qui met en dualit\'e les espaces ${\got u}_p$ et 
${\got u}_{p-1}$. ${\mathcal F}_p$ se prolonge en un morphisme d'alg\`ebres 
${\mathcal F}_p : {\mathbb D}({\mathbb X}_{p-1}) \lra {\mathbb D}({\mathbb X}_p)$ et on pose,
enfin : $\phi_{p-1} = {\mathcal F}_p \circ \psi_{p-1}$.
On obtient ainsi un morphisme d'alg\`ebres $\phi_{p-1} : U({\got p}_{p-1}) 
\lra {\mathbb D}({\mathbb X}_p)$ et il est clair que :
$$\ker \phi_{p-1} = I_{\pi,p-1,p}, \ \ker \phi_p = I_ {\pi,p,p}$$
Il suffit, pour conclure, de v\'erifier que les restrictions de $\phi_p$ et 
$\phi_{p-1}$ \`a $P_{p-1,p}$ sont \'egales. Compte-tenu des d\'efinitions et
 par analogie avec le cas pr\'ec\'edent, il suffit de prouver que :
$$\phi_p(Z) = \phi_{p-1} (Z), \forall Z \in <X_{\alpha_{p+1}}, 
X_{-\alpha_{p+1}}, H_{\alpha_{p+1}}, X_{\alpha_p}, X_{-\alpha_p}, 
H_{\alpha_p}>$$
Or, ceci se v\'erifie sans difficult\'es, en tenant compte des choix faits dans les 
d\'efinitions de $E_p^\infty, E_{p-1}^\infty$ et de l'espace de 
r\'ealisation de $\rho_{-i}$.

\smallskip

{\bf Th\'eor\`eme 6.8 :} {\it Soit $(k,\var) \in I_n$. Soit $\pi_k
\in {\mathcal R}_k$. 
Alors :
$$GKdim \pi_k = \dim O_{k,\var}$$
Ainsi, la repr\'esentation $\pi_k$ est $GK$-associ\'ee \`a l'orbite $O_{k,\var}$.}

\smallskip

{\bf Preuve :} D'apr\`es la proposition 6.5, on dispose des morphismes 
d'alg\`ebres $\phi_i : U({\got p}_i) \lra {\mathbb D}({\mathbb X}_k),      
\phi_j : U({\got p}_j) \lra {\mathbb D}({\mathbb X}_k), (i,j) \in \{k-1,k\}$ ou $(i,j) 
\in \{k,k+1\}$ tels que~:
$$(\phi_i)_{\mid U({\got p}_{ij})} = (\phi_j)_{\mid U({\got p}_{ij})}$$
D'apr\`es la proposition 6.1 et le corollaire 6.2, l'alg\`ebre ${\mathbb D}(
{\mathbb X}_k)$ est solution du probl\`eme universel pour la somme amalgam\'ee 
des $U({\got p}_i), i \in \{k-1,k\}$ ou $i \in \{k,k+1\}$ et il existe donc 
un morphisme d'alg\`ebres  $\phi : U({\got g}) \lra {\mathbb D}({\mathbb X}_k)$ 
tel que :
$$\phi_{\mid U({\got p}_i)} = \phi_i, \ i \in \{k-1,k\} \ {\rm ou} \ i \in 
\{k,k+1\}$$
Soit $E_k$ l'espace de la repr\'esentation $\pi_k$.
En utilisant (21), (22), (35), (37), (41) et (42), on peut \'ecrire :
$$\forall f \in E_k^\infty, \forall a \in U({\got p}_i), i \in \{k-1,k\} \
 {\rm ou} \ i \in \{k,k+1\}, 
\ \pi_k^\infty (a).f = \phi_i(a).f$$
Comme $\got g$ est engendr\'ee par les deux paraboliques ${\got p}_i, i \in 
\{k-1,k\}$ ou $i \in \{k,k+1\}$, on en d\'eduit que~:
$$\forall f \in E_k^\infty, \ \forall a \in U({\got g}), 
\pi_k^\infty (a).f = \phi(a).f$$
Il s'en suit que $\phi$ induit un morphisme injectif $\dis{\overline{\phi} :
(U({\got g}) / \ker \pi_k^\infty ) \lra {\mathbb D}({\mathbb X}_k)}$.

Ceci implique que : 
$$GKdim \pi_k \leq GKdim \ {\mathbb D}({\mathbb X}_k)$$
Or, selon des r\'esultats classiques (voir, par exemple, [SM]), on sait que :
$$GKdim \ {\mathbb D}({\mathbb X}_k) = 2 \dim {\mathbb X}_k = \dim O_{k,\var}$$       
Ceci implique l'in\'egalit\'e : $\dis{GKdim \pi_k \leq  
\dim O_{k,\var}}$. 
L'\'egalit\'e souhait\'ee est alors cons\'e\-quence de (20).
\qed

\section{\bf Repr\'esentations $GK$-associ\'ees et repr\'esentations
associ\'ees.}
   
La m\'ethode des orbites consiste \`a ``associer'', selon un sens
qui sera rappel\'e ult\'e\-rieu\-rement, le dual unitaire de $G$ et
l'ensemble des $G$-orbites nilpotentes coadjointes. Le but de ce
paragraphe est donc de montrer que, pour tout entier $k$, la
repr\'esentation $\pi_k$ est ``associ\'ee'' \`a $O_{k,\var}$. 

\subsection{} Pla\c cons-nous \`a nouveau dans le cadre
g\'en\'eral suivant :

Soit $G$ un groupe r\'eel simple connexe et simplement connexe d'alg\`ebre de Lie $\got g$,
${\got g}_{\mathbb C}$ la complexifi\'ee de $\got g$, $G_{\mathbb C}$ un
groupe complexe simplement connexe d'alg\`ebre de Lie ${\got g}_{\mathbb
C}, U({\got g})$
l'alg\`ebre enveloppante de ${\got g}_{\mathbb C}$ et $S({\got g})$ son
alg\`ebre sym\'etrique. On suppose choisie
une sous-alg\`ebre de Borel $\got b$ de $\got g$ et on note ${\got
b}_{\mathbb C}$ sa complexifi\'ee. Soit $B$ le sous-groupe de Borel de
$G$, d'alg\`ebre de Lie $\got b$. Consid\'erons maintenant une
repr\'esentation unitaire irr\'eductible $\pi$ de $G$, $\pi_B$ sa
restriction \`a $B$. On note $I_{\pi,G}$
l'annulateur infinit\'esimal de $\pi$ dans $U({\got g})$, $I_{\pi,B} =
I_{\pi,G} \cap U({\got b})$ l'annulateur infinit\'esimal de $\pi_B$.
Soit, enfin, $Gr : U({\got g}) \lra S({\got g})$ l'application
canonique usuelle qui envoie $U({\got g})$ sur son gradu\'e $S({\got
g})$. On introduit de m\^eme l'application $Gr : U({\got b}) \lra
S({\got b})$.

On pose : $J_{\pi,G} = Gr(I_{\pi,G}), J_{\pi,B} = Gr(I_{\pi,B}) =
J_{\pi,G} \cap S({\got b})$. Soit $V(J_{\pi,G})$ la vari\'et\'e des
z\'eros dans ${\got g}^*_{\mathbb C}$ de $J_{\pi,G}$.
Si $I_{\pi,B}$ est un id\'eal primitif de $U({\got b})$, on sait
d\'efinir un id\'eal de $S({\got b})$,
not\'e $\widetilde{J}_{\pi,B}$, image de $I_{\pi,B}$ par l'inverse de
l'application 
de Dixmier. Nous reviendrons pr\'ecis\'ement sur ce point dans 7.2.
 
La notion de repr\'esentation associ\'ee \`a une $G$-orbite
r\'eelle est, selon la terminologie usuelle de la m\'ethode des
orbites, celle qui est rappel\'ee dans la d\'efinition suivante :

\smallskip

{\bf D\'efinition 7.1:} {\it Soit $O$ une $G$-orbite coadjointe
dans ${\got g}^*$ et $O_{\mathbb C}$ sa complexifi\'ee. $\pi$ est dite
associ\'ee \`a $O$ si:
$$V(J_{\pi,G}) = \overline{O_{\mathbb C}}$$ }

\smallskip

On rappelle enfin que, d'apr\`es \cite{BB}, il existe une $G_{\mathbb C}$-orbite
adjointe nilpotente $O_{\pi}$ dans ${\got g}_{\mathbb C}$ telle que :
$$V(J_{\pi,G}) = \overline{O_\pi}$$
Nous allons commencer par d\'emontrer le r\'esultat suivant~:

\smallskip

{\bf Th\'eor\`eme 7.1 :} {\it On suppose que la repr\'esentation
$\pi$ satisfait aux trois hypoth\`eses suivantes~:

- (1) La repr\'esentation $\pi_B$ est irr\'eductible.

- (2) L'id\'eal $\widetilde{J}_{\pi,B}$ est un id\'eal gradu\'e de
$S({\got b})$.

- (3) $GKdim \pi = GKdim \pi_B = GKdim (S({\got b}) / \widetilde{J}_{\pi,B})$.

Alors, la $G_{\mathbb C}$-orbite $O_{\pi}$ contient une $B_{\mathbb
C}$-orbite ouverte.}

\subsection{} La d\'emonstration du th\'eor\`eme 7.1. repose
fortement sur un r\'esultat de J.Y.Char\-bon\-nel, donn\'e dans
\cite{CH}, concernant les id\'eaux
primitifs d'une alg\`ebre de Lie compl\`etement r\'e\-solu\-ble. Pour
simplifier, et puisque nous appliquerons ce qui va suivre \`a la
sous-alg\`ebre de Borel introduite pr\'ecedemment, nous noterons $\got
b$ une telle alg\`ebre de Lie et nous la supposerons alg\'ebrique de
radical unipotent $\got n$. 

Soit $A({\got b})$ l'alg\`ebre des
op\'erateurs diff\'erentiels sur $\got b$ \`a coefficients polyn\^omiaux 
, $\widehat{S}({\got b}^*)$ l'alg\`ebre des s\'eries formelles sur ${\got
b}$,  $\widehat{A}({\got b})$
l'alg\`ebre des op\'erateurs diff\'erentiels sur $\got b$ \`a
coefficients s\'eries formelles, $E_{\got b}$ le sous-${\mathbb
Q}$-espace vectoriel de ${\got b}^*$ engendr\'e par les poids de la
repr\'esentation adjointe de $\got b$, $\widehat{S}(E_{\got b})$ le
sous-anneau ferm\'e de $\widehat{A}({\got b})$ engendr\'e par
$E_{\got b}$. Soit, enfin, $\widehat{P}({\got b})$ la sous-alg\`ebre
de $\widehat{A}({\got b})$ engendr\'ee par $A({\got b})$ et
$\widehat{S}(E_{\got b})$.

On utilise pour $\widehat{P}({\got b})$ la filtration 
induite par la filtration naturelle d\'efinie sur $\widehat{A}({\got
b})$ et on consid\`ere le gradu\'e associ\'e qui s'identifie \`a
$\widehat{S}(E_{\got b}).S({\got b}^*) \otimes S({\got b})$. Si $L$
est un id\'eal \`a gauche de $\widehat{P}({\got b})$, on note $Gr(L)$
l'id\'eal correspondant dans $\widehat{S}(E_{\got b}).S({\got b}^*)
\otimes S({\got b})$. 

Dans \cite{DI2}, J.Dixmier introduit les op\'erateurs $L_{\got b},
R_{\got b}$ et $W_{\got b}$ d\'efinis de la mani\`ere suivante;
On consid\`ere les deux s\'eries enti\`eres :
$$\frac{T}{1 - \exp (-T)} = \sum b_r T^r, \ \ \frac{T}{\exp (T) - 1} =
\sum c_r T^r$$ 
On v\'erifie que : $\dis{b_0 = c_0 = 1, b_1 = \frac{1}{2}, c_1 = -\frac{1}
{2}, b_r = c_r, \forall r \geq 2, b_{2r+1} = 0}$.
Pour tout $x \in {\got b}$, on pose~:
$$\begin{array}{rl}
\forall y \in {\got b}, \ L_{\got b}(x)(y) &= \sum_{r \geq 0}
b_r (ad\ y)^r.x \\
\ R_{\got b}(x)(y) &= \sum_{r \geq 0}
c_r (ad\ y)^r.x \\
W_{\got b}(x)(y) &= L_{\got b}(x)(y) - R_{\got b}(x)(y) = [y,x] 
\end{array}
$$

L'application $L_{\got b}$ est un homomorphisme d'alg\`ebres de Lie de
$\got b$ dans l'alg\`ebre de Lie sous-jacente \`a $\widehat{A}({\got
b})$. On peut prolonger $L_{\got b}$ de mani\`ere canonique en un
homomorphisme de $U({\got b})$ dans $\widehat{A}({\got b})$.

Soit $(e_1, \dots,e_n)$ une base de $\got b$, $(e^*_1,
\dots,e^*_n)$ la base duale dans ${\got b}^*$. J.Dixmier, dans
\cite{DI2}, lemme 7.3,
donne une expression des champs de vecteurs $L_{\got b}$ \`a partir
des bases pr\'ec\'edentes :

\smallskip

{\bf Lemme 7.2 :} {\it Soit $\dis{u =
\sum_{|\alpha| \leq p} \lambda_\alpha e_1^{\alpha_1}\dots
e_n^{\alpha_n}}$ un \'el\'ement de $U_p({\got b})$ o\`u les
$\lambda_\alpha $ sont des nombres complexes. Alors, on a :
$$L_{\got b}(u) = \sum_{|\alpha| = p} \lambda_\alpha
e_1^{\alpha_1}\dots e_n^{\alpha_n} + \sum_{|\beta| < p} \psi_\beta
e_1^{\beta_1}\dots e_n^{\beta_n} + \sum_{|\gamma| < p, 1 \leq i \leq n}
w_{\gamma,i}e_1^{\gamma_1}\dots e_n^{\gamma_n}W_{\got b}(e_i)$$
o\^u les $\psi_\beta, w_{\gamma,i}$ sont des \'el\'ements de
$\widehat{S}({\got b}^*)$.} 

\smallskip

Dans \cite{DI1}, J.Dixmier introduit une application dite
``application de Dixmier `` et not\'ee $Dix_{\got b}$ qui, \`a chaque
id\'eal premier $\got g$-invariant $J$ de $S({\got b})$, associe un
id\'eal premier $Dix_{\got
b}(J)$ de $U({\got b})$. En fait, $Dix_{\got g}(J)$ se r\'ealise comme l' annulateur
d'une repr\'esentation induite de $\got b$ \`a partir
d'une polarisation d'un \'el\'ement $f$ de ${\got b}^*$ et $J$ n'est autre
que l'id\'eal de l'orbite $B_{\mathbb C}.f$. Cette orbite sera appel\'ee
``l'orbite de Dixmier'' de $Dix_{\got g}(J)$.

L'application $Dix_{\got b}$ est une bijection de l'ensemble des
id\'eaux premiers $\got g$-invariants de $S({\got b})$ sur l'ensemble
des id\'eaux premiers de $U({\got b})$. On notera $\beta_{\got b}$ son
application inverse.  

Si $I$ est un id\'eal de $U({\got b})$, on note enfin $\psi_{\got
b}(I)$ l'id\'eal \`a gauche de $\widehat{P}({\got b})$ engendr\'e par
$L_{\got b}(I)$ et $W_{\got b}({\got n})$. On a alors le r\'esultat
suivant, donn\'e dans \cite{CH}, corollaire 5.8 :

\smallskip

{\bf Th\'eor\`eme 7.3 :} {\it Soit $I$ un
id\'eal primitif de $U({\got b})$ et $J = \beta_{\got b}(I)$. Alors,
on a l'inclusion suivante :
$$\sqrt{Gr(J)} \subset Gr(\psi_{\got b}(I))$$}

On note $\phi$  le morphisme surjectif canonique d'alg\`ebres de $\widehat{A}({\got
b})$ sur $S({\got b})$ d\'efini de la mani\`ere suivante : 
Soit $\dis{u \in \widehat{A}({\got b}), u = \sum_{|\alpha| \leq p}
\varphi_\alpha e_1^{\alpha_1}\dots e_n^{\alpha_n}}$, avec
$\varphi_\alpha \in \widehat{S}({\got g}^*)$, pour tout $n$-uplet
$\alpha$. Alors, 
$$\phi(u) = \sum_{|\alpha| \leq p} \varphi_\alpha(0)
e_1^{\alpha_1}\dots e_n^{\alpha_n}$$
Comme $W_{\got b}(e_i) = e_i^*.e_i$, il s'en suit que : $\forall i, 1
\leq i \leq n, \phi(W_{\got b}(e_i)) = 0$.

D'autre part, soit $\dis{u = \sum_{|\alpha| \leq p} \lambda_\alpha
e_1^{\alpha_1} \dots e_n^{\alpha_n} \in U_p({\got b})}$. A l'aide du
lemme 7.2 on peut donc \'ecrire :
$$\phi(L_{\got b}(u)) = \sum_{|\alpha| = p} \lambda_\alpha
e_1^{\alpha_1}\dots e_n^{\alpha_n} + \sum_{|\beta| < p} \psi_\beta(0)
e_1^{\alpha_1}\dots e_n^{\alpha_n}$$
De ceci on d\'eduit que, si $Gr(\phi)$ d\'esigne l'application
gradu\'ee correspondante, alors on a :
$$\forall u \in U({\got b}), Gr(\phi)((Gr(L_{\got b}(u))) = Gr(u)$$
Ceci implique que, pour tout  id\'eal $I$ de $U({\got b})$ :
$$Gr(\phi)(Gr(\psi_{\got b}(I))) \subset Gr(I)$$
Soit $I$ un id\'eal primitif de $U({\got b}), J = \beta_{\got
b}(I)$. On sait, d'apr\`es le th\'eor\`eme 7.3, que $\sqrt{Gr(J)}
\subset Gr(\psi_{\got b}(I))$. Comme, d'autre part, $\sqrt{Gr(J)}$ est
un id\'eal de $S({\got b})$, il s'en suit que $Gr(\phi)(\sqrt{Gr(J)})
= \sqrt{Gr(J)}$. D'o\`u :
$$\sqrt{Gr(J)} \subset Gr(I)$$
Ainsi, on obtient :

\smallskip

{\bf Corollaire 7.4 :} {\it Soit $I$ un id\'eal primitif de
$U({\got b})$ et $J = \beta_{\got b}(I)$. Alors, on a :
$$\sqrt{Gr(J)} \subset Gr(I)$$}

\subsection{\bf Preuve du th\'eor\`eme 7.1.} Supposons donn\'ee une
repr\'esentation unitaire irr\'e\-duc\-tible $\pi$ de $G$ v\'erifiant les
hypoth\`eses du th\'eor\`eme 7.1. La repr\'esentation $\pi_B$ est
irr\'eductible, l'id\'eal $I_{\pi,B}$ est donc primitif et comme
, d'apr\`es l'hypoth\`ese (2) du th\'eor\`eme 7.1, $\widetilde{J}_{\pi,B}$ est gradu\'e dans $S({\got b})$, on a,
d'apr\`es le corollaire 7.4 :
$$\sqrt{\widetilde{J}_{\pi,B}} \subset Gr(I_{\pi,B})$$
Soit $O_{\pi,B}$ la $B_{\mathbb C}$-orbite de Dixmier de l'id\'eal
primitif $I_{\pi,B}$. On sait que
$\widetilde{J}_{\pi,B}$ est l'id\'eal de $O_{\pi,B}$. Comme
$\overline{O_{\pi,B}}$ est une sous-vari\'et\'e irr\'eductible de
${\got b}^*_{\mathbb C}$ il s'en suit que $\widetilde{J}_{\pi,B}$ est un
id\'eal premier, ce qui implique :
$$\widetilde{J}_{\pi,B} \subset Gr(I_{\pi,B})$$
D'apr\`es l'hypoth\`ese (3) du th\'eor\`eme 7.1, on sait que les
id\'eaux $\widetilde{J}_{\pi,B}$ et $Gr(I_{\pi,B})$ ont m\^eme
dimension de Gelfand-Kirillov ou encore que~:
$$GKdim(S({\got b}) /
\widetilde{J}_{\pi,B}) = GKdim(S({\got b}) / Gr(I_{\pi,B}))$$
 Comme
$\widetilde{J}_{\pi,B}$ est premier, l'inclusion pr\'ec\'edente et le
corollaire 3.6 de \cite{BO-KR} permettent alors d'affirmer que :
\begin{equation}
\widetilde{J}_{\pi,B} = Gr(I_{\pi,B}) = J_{\pi,B} 
\end{equation} 

Soit $\mu : {\got g}^*_{\mathbb C} \lra {\got b}^*_{\mathbb C}$
l'application ``restriction'' et $\mu^* : S({\got b}) \lra S({\got
g})$ le comorphisme de $\mu$, identifiant $S({\got b})$ \`a une
sous-alg\`ebre de $S({\got g})$. 

Via la forme de Killing, on peut donc
d\'efinir un morphisme dominant, not\'e encore $\mu$,  de $O_\pi$ sur
la sous-vari\'et\'e ferm\'ee $\overline{\mu(O_\pi)}$ de ${\got
b}^*_{\mathbb C}$. On d\'esigne  respectivement par $R(\overline{O_\pi}), R(\overline{\mu(O_\pi)})$ et
$R(\overline{O_{\pi,B}})$ les alg\`ebres de fonctions sur les
vari\'et\'es affines correspondantes. 

$J_{\pi,G}$ est l'id\'eal, dans $S({\got g})$, de la vari\'et\'e
$\overline{O_\pi}$  et il est facile d'en d\'eduire que, moyennant
l'identification pr\'ec\'edente, $J_{\pi,B}$ est l'id\'eal, dans
$S({\got b})$,  de la vari\'et\'e $\overline{\mu(O_\pi)}$.

Il existe un morphisme canonique
surjectif de $S({\got b})$ sur $R(\overline{\mu(O_\pi)})$ qui, \`a chaque
polyn\^ome $P$ sur ${\got g}^*_{\mathbb C}$, fait correspondre la restriction de
$P$ \`a $\overline{\mu(O_\pi)}$, consid\'er\'ee comme sous-vari\'et\'e de
${\got b}^*_{\mathbb C}$. De ce qui pr\'ec\`ede on d\'eduit un isomorphisme d'alg\`ebres 
 $F :S({\got b}) / J_{\pi,B} \lra R(\overline{\mu(O_\pi)})$.

Il existe aussi un morphisme surjectif canonique de $S({\got b})$ sur
$R(\overline{O_{\pi,B}})$ qui, \`a tout polyn\^ome $P$ sur ${\got
b}_{\mathbb C}^*$, fait correspondre la restriction de $P$ \`a
$\overline{O_{\pi,B}}$ et ce morphisme  
induit un  isomorphisme
d'alg\`ebres  de 
$S({\got b}) / \widetilde{J}_{\pi,B}$ sur
$R(\overline{O_{\pi,B}})$. Comme $\widetilde{J_{\pi,B}} = J_{\pi,B}$, 
On  d\'eduit de tout ceci l'existence d'un isomorphisme
d'alg\`ebres  $\sigma^* : R(\overline{\mu(O_\pi)}) \lra
R(\overline{O_{\pi,B}})$. On constate, de plus, que cet isomorphisme
est $B_{\mathbb C}$-\'equivariant.
Suivant des r\'esultats bien connus de g\'eom\'etrie alg\'ebrique,
(voir, par exemple, \cite{HA}, corollaire 3.7), on
 sait qu'il existe un isomorphisme de vari\'et\'es $\sigma :
\overline{O_{\pi,B}} \lra \overline{\mu(O_\pi)}$ dont $\sigma^*$ est
le comorphisme. 

Il est facile de v\'erifier, en outre, que cet
isomorphisme est $B_{\mathbb C}$-\'equivariant. En effet, soit $u \in 
\overline{O_{\pi,B}}, b \in B_{\mathbb C}$ et supposons que $\sigma (b.u)
\not= b.\sigma(u)$. Il existe alors $f \in R(\overline{\mu(O_\pi)})$
tel que : $f(\sigma(b.u)) \not= f(b.\sigma(u))$, ce qui revient \`a
dire que $b^{-1}.\sigma^*(f)(u) \not= \sigma^*(b^{-1}.f)(u)$. Ceci
contredit alors la $B_{\mathbb C}$-\'equivariance de $\sigma^*$. 
On obtient ainsi un morphisme de vari\'et\'es $\theta = \sigma^{-1}
\circ \mu : O_{\pi} \lra \overline{O_{\pi,B}}$ qui est dominant et
$B_{\mathbb C}$-\'equivariant.

De l'hypoth\`ese (3) du th\'eor\`eme 7.1., on d\'eduit que : $\dim
O_{\pi} = \dim O_{\pi,B}$. D'autre part, l'orbite $O_{\pi}$ est une r\'eunion de
$B_{\mathbb C}$-orbites. Supposons que chaque $B_{\mathbb C}$-orbite dans
$O_\pi$ soit de dimension strictement inf\'erieure \`a $\dim
O_\pi$. Soit $O$ une telle orbite . Comme $\theta$ est $B_{\mathbb
C}$-\'equivariant, il s'en suit que $\theta(O)$ est une $B_{\mathbb C}$
-orbite dans $\overline{O_{\pi,B}}$, de dimension strictement
inf\'erieure \`a $\dim O_{\pi,B}$. D'o\`u : $\theta(O) \subset
\overline{O_{\pi,B}} \backslash O_{\pi,B}$. Ceci est vrai pour toute
$B_{\mathbb C}$-orbite $O$ dans $O_\pi$. On a, ainsi :
$$\theta(O_{\pi}) \subset \overline{O_{\pi,B}} \backslash
O_{\pi,B}$$
Comme $O_{\pi,B}$ est ouvert dans $\overline{O_{\pi,B}}$, cela
contredit le fait que le morphisme $\theta$ est dominant. Ainsi,
$O_{\pi}$ poss\`ede une $B_{\mathbb C}$-orbite de dimension $\dim O_\pi$,
ce qui finit de d\'emontrer le th\'eor\`eme 7.1.

\subsection{}On va maintenant appliquer le th\'eor\`eme 7.1 au cas
des repr\'esentations $\pi_k$ construites
pr\'ec\'edemment. On reprend les notations des paragraphes 3,4,5 et 6.

Soit $(k,\var) \in I_n ,\pi_k \in {\mathcal R}_k, I_{\pi,k}$
l'annulateur infinit\'esimal de $\pi_k$ dans $U({\got g})$ et
$J_{\pi,k} = Gr(I_{\pi,k})$ l'id\'eal correspondant dans $S({\got
g})$. On d\'esigne par $O_{k,{\mathbb C}}$ la $G_{\mathbb C}$-orbite telle
que :
$$V(J_{\pi,k}) = \overline{O_{k,{\mathbb C}}}$$

- D'apr\`es la proposition 5.4, la restriction $\pi_{k,B}$ de
$\pi_k$ \`a $B$ est irr\'eductible.

- Soit $b_{k,\var} = {\mathcal K}(X_{k,\var}, .)$ la restriction \`a
${\got b}^*$ de la forme lin\'eaire associ\'ee \`a $X_{k,\var}$ par la
forme de Killing sur $\got g$. $I_{\pi,k,B}$ est l'annulateur
infinit\'esimal de $\pi_{k,B}$ et il r\'esulte de la construction de
$\pi_{k,B}$ que cet id\'eal est exactement celui qui correspond \`a
l'orbite $O_{\pi,k,B} = B_{\mathbb C}.b_{k,\var}$ par la correspondance
de Dixmier (voir \cite{DI1}, chapitre 6). Soit $\widetilde{J}_{\pi,k,B}$ l'id\'eal de
$O_{\pi,k,B}$.
Comme $\dim O_{\pi,k,B} = \dim O_{k,\var}$, il s'en suit que :
$$GKdim \pi_k = GKdim \ (S({\got b}) / \widetilde{J}_{\pi,k,B})$$

- En reprenant les arguments de 6.2. et en utilisant notamment la
proposition 6.3. on montre \'egalement que :
$$GKdim \pi_{k,B} = \dim O_{k,\var} = GKdim \pi_k$$   

- On sait, enfin, que $b_{k,\var}$ est d\'efinie par un \'el\'ement
nilpotent de $\got g$. Ceci implique facilement  que l'orbite
$O_{\pi,k,B}$ est un c\^one et, donc, que son id\'eal
$\widetilde{J}_{\pi,k,B}$ est gradu\'e.

Ainsi, les trois hypoth\`eses du th\'eor\`eme 7.1. sont
satisfaites et on obtient : 

\smallskip

{\bf Proposition 7.5 :} {\it Soit $(k,\var) \in I_n, \pi_k$ un
\'el\'ement de ${\mathcal R}_k$. Alors, la $G_{\mathbb C}$-orbite associ\'ee
\`a $\pi_k$ contient une $B_{\mathbb C}$-orbite dense.}

\smallskip

On sait maintenant que, dans $sl_n({\mathbb C})$, il n'existe qu'une
seule orbite nilpotente sph\'e\-ri\-que de dimension donn\'ee.  La
proposition 7.5 et le th\'eor\`eme 6.8 impliquent  finalement le
r\'esultat suivant~:

\smallskip

{\bf Th\'eor\`eme 7.6 :} {\it Soit $(k,\var) \in I_n, \pi_k \in
{\mathcal R}_k$. Alors, la repr\'esentation $\pi_k$ est associ\'ee \`a
l'orbite $O_{k,\var}$.}

\affiliationone{Herv\'e Sabourin \\
UMR 6086 CNRS\\
D\'epartement de Math\'ematiques\\
Universit\'e de Poitiers\\
Boulevard Marie et Pierre Curie\\
T\'el\'eport 2 - BP 30179\\
86962 Futuroscope Chasseneuil cedex\\
France\\
sabourin@math.univ-poitiers.fr}

\end{document}